\documentclass{article} 

\usepackage[utf8]{inputenc} 
\usepackage[T1]{fontenc}    
\usepackage{hyperref}       
\usepackage{url}            
\usepackage{booktabs}       
\usepackage{amsfonts}       
\usepackage{nicefrac}       
\usepackage{microtype}      
\usepackage{color,xcolor}         

\usepackage{fullpage}

\usepackage{natbib} 
\usepackage{mathtools}
\newcommand{\eqdef}{\coloneqq}

\usepackage{algorithm}
\usepackage{algorithmic}

\usepackage{amsmath, amsfonts,amssymb,amsthm}
\usepackage{cite}
\DeclareMathOperator*{\argmin}{arg\,min}

\newcommand{\sqn}[1]{\left\| #1 \right\|^2}
\newcommand{\Exp}[1]{\mathbb{E}\!\left[ #1 \right]} 
\newcommand{\Expb}[2]{\mathbb{E}\!\left[#1\ \vert\ #2\right]}
\newcommand{\Expc}[1]{\Expb{#1}{\mathcal{F}^t}}
\newcommand{\Exps}[2]{\mathbb{E}_{#1}\left[ #2 \right]} 

\theoremstyle{plain}
\newtheorem{theorem}{Theorem}[section]

\newtheorem{corollary}[theorem]{Corollary}
\theoremstyle{definition}

\newtheorem{assumption}[theorem]{Assumption}
\theoremstyle{remark}

\usepackage{xspace}
\usepackage{scalefnt}
\definecolor{mydarkredd}{rgb}{0.7,0.0,0.0}
\newcommand{\algn}[1]{{\sf\color{mydarkredd}\scalefont{0.955}{#1}}\xspace}
\newcommand{\algno}{\algn{SMPM}}
\newcommand{\algnod}{\algn{FedSMPM}}

\definecolor{lgray}{rgb}{0.95,0.95,0.95}
\definecolor{yel}{rgb}{1,0.98,0.92}
\definecolor{mydarkbluee}{rgb}{0,0.2,0.9}
\newcommand\bl[1]{{\color{mydarkbluee}#1}}
\newcommand\ff{\bl{f}}

\definecolor{mydarkred}{rgb}{0.8,0.0,0.0}
\newcommand\dr[1]{{\color{mydarkred}#1}}
\newcommand\g{\dg{g}}

\definecolor{mydarkgreen}{rgb}{0,0.55,0}
\newcommand\dg[1]{{\color{mydarkgreen}#1}}
\newcommand\h{\dr{h}}

\newcommand\hi{\dr{h_i}}

\newcommand{\vv}{u}

\newcommand{\RR}{\mathbb{R}}
\newcommand{\bW}{\mathbf{W}}

\author{Laurent Condat$^{1,2}$ \qquad Elnur Gasanov$^{1}$ \qquad    Peter Richt\'{a}rik$^{1,2}$\\[3mm]
$^1$Computer Science Program, CEMSE Division,\\ King Abdullah University of Science and Technology (KAUST)\\ Thuwal, 23955-6900, Kingdom of Saudi Arabia\\
 $^2$SDAIA-KAUST Center of Excellence in Data Science and \\Artificial Intelligence 
(SDAIA-KAUST AI)\\[3mm]
Contact: see \texttt{https://lcondat.github.io/} 
}

\title{\textbf{The Stochastic Multi-Proximal Method \\for Nonsmooth Optimization}}

\date{January 31, 2025}

\begin{document}

\maketitle

\begin{abstract}
Stochastic gradient descent type methods are ubiquitous in machine learning, but they are only applicable to the optimization of differentiable functions. Proximal algorithms are more general and applicable to nonsmooth functions. We propose a new stochastic and variance-reduced algorithm, the Stochastic Multi-Proximal Method (\algno), in which the proximity operators of a (possibly empty) random subset of functions are called at every iteration, according to an arbitrary sampling distribution. Several existing algorithms, including  \algn{Point-SAGA} (2016), \algn{Proxskip}  (2022) and \algn{RandProx-Minibatch} (2023) are recovered as particular cases. We derive linear convergence results in presence of strong convexity and smoothness or similarity of the functions. We prove convergence in the general convex case and accelerated $\mathcal{O}(1/t^2)$  convergence with varying stepsizes in presence of strong convexity solely. Our results are new even for the above special cases. Moreover, we show an application to distributed optimization with compressed communication, outperforming existing methods.
\end{abstract}

{
\renewcommand\baselinestretch{0}
\tableofcontents
\renewcommand\baselinestretch{1}
}

\section{Introduction}

Optimization problems arise in virtually every field and are pervasive in machine learning (ML) and artificial intelligence (AI) \citep{pal09, sra11, bac12, cev14, pol15, bub15, glo16, cha16, sta16}, notably for training large models. Because of the high dimension of the variables to optimize and the huge size of the datasets to be scanned during the iterative optimization process, stochastic first-order methods are well suited, and Stochastic Gradient Descent (\algn{SGD}) has become the cornerstone of modern ML and AI applications \citep{RobbinsMonro:1951,Nemirovski-Juditsky-Lan-Shapiro-2009,Bottou2012,kha20}. But \algn{SGD} is only applicable to differentiable functions, for which the gradient is defined and regular enough. On the other hand, proximal algorithms \citep{par14,ryu16,con23}, which use the proximity operator (prox) of functions, are more general and applicable to nonsmooth and even infinite-valued functions. We recall that the prox of a function $\phi$ is the operator
 $
 \mathrm{prox}_{\phi}: x \mapsto \argmin_{y}\big(\phi(y) + \frac{1}{2} \sqn{y-x}\big)
  $  \citep{bau17}. 
  This operator has a closed form for many functions of practical interest \citep{par14,pus17,ghe18}, see also the website \href{http://proximity-operator.net}{proximity-operator.net}. But even if a prox can be computed, it might be expensive to do so. So, reducing the number of prox calls needed to achieve the desired accuracy in iterative algorithms is crucial.

We propose a new stochastic algorithm for convex nonsmooth optimization, the \textbf{Stochastic Multi-Proximal Method} (\algno), 
in which the proxs of a (possibly empty) random subset of functions are called at every iteration, according to an \textbf{arbitrary sampling} distribution. \algno is variance-reduced \citep{han19,gor202,gow20a}; that is, it converges to an exact solution, thanks to its primal--dual nature. Several existing algorithms, including \algn{Point-SAGA}
 \citep{def16}, 
 the \algn{Stochastic Decoupling Method} \citep{mis19}, 
 \algn{Proxskip} \citep{mis22} and \algn{RandProx-Minibatch} \citep{con23rp} are recovered as particular cases.
We derive linear convergence results in presence of strong convexity, and smoothness \textbf{or} similarity. It is remarkable to obtain linear convergence in settings where the functions are not smooth. We also derive accelerated $\mathcal{O}(1/t^2)$ rates with varying stepsizes for strongly convex problems, without other assumptions, and we prove convergence in the general convex case. Our results are new even for the aforementioned particular cases. Moreover, we show an application to distributed optimization with compressed communication, outperforming existing methods.

\subsection{Problem Formulation}

Given $n\geq 1$, we consider the generic convex optimization problem
\begin{equation}
\mathrm{Find} \ x^\star \in \argmin_{x\in\mathcal{X}}  \bigg( \ff(x) + \g(x)+\frac{1}{n}\sum_{i=1}^n \hi(x)\bigg)\label{eqpbm}
 \end{equation}
 involving proper  closed convex functions $\ff:\mathcal{X}\rightarrow \mathbb{R}$, $\g:\mathcal{X}\rightarrow \mathbb{R}\cup\{+\infty\}$,
 $\hi:\mathcal{X}\rightarrow \mathbb{R}\cup\{+\infty\}$,  $i\in[n]\eqdef\{1, \ldots,n\}$, 
  defined on a finite-dimensional real Hilbert space $\mathcal{X}$ \citep{bau17}, for instance $\mathbb{R}^d$ for a large dimension $d\geq 1$. 
 The generic template \eqref{eqpbm} captures a large variety of problems with a wide range of applications. Notably, the constraint that $x$ belongs to some set $\Omega$ can be included as a function in the problem by letting $\g$ or one $\hi$ be the indicator function $\imath_{\Omega} = \{0$ if $x\in \Omega$, $+\infty$ otherwise$\}$.  
Then   $\imath_{\Omega}$ is proper closed and convex if $\Omega$ is nonempty closed and convex, and $\mathrm{prox}_{\imath_{\Omega}}$ is the projection onto $\Omega$ \citep{bau17}. The following general assumption is supposed to hold throughout the paper for the problem \eqref{eqpbm} to be well posed: 
 $\ff$ is $L_{\ff}$-smooth\footnote{A function $\phi:\mathcal{X}\rightarrow \mathbb{R}$ is said to be $L$-smooth if it is differentiable and its gradient is $L$-Lipschitz continuous; that is, for every $x\in\mathcal{X}$ and $y\in\mathcal{X}$, $\|\nabla \phi(x)-\nabla \phi(y)\|\leq L \|x-y\|$.}, for some $L_{\ff}\geq 0$. Also, 
there exists $x^\star\in\mathcal{X}$ and $u_i^\star \in \partial\hi(x^\star)$, $\forall i\in[n]$, such that
\begin{equation}
\textstyle
0\in \nabla \ff(x^\star)+\partial \g (x^\star)+ \frac{1}{n}\sum_{i=1}^n u_i^\star,\label{eqincl}
\end{equation}
where $\partial\cdot$ denotes the subdifferential \citep{bau17}.
Such an $x^\star$ is a solution to \eqref{eqpbm}.
 We denote by $\mu_{\ff}\geq 0$, $\mu_\g\geq 0$, $\mu_{\hi}\geq 0$ the
strong convexity constants of $\ff$, $\g$, and $\hi$,
  respectively ($\phi$ is said to be $\mu$-strongly convex if $\phi-\frac{\mu}{2}\sqn{\cdot}$ is convex). 
 For every $i\in[n]$, we let $L_\hi >0$ be such that $\hi$ is $L_\hi$-smooth; if  $\hi$ is not smooth, we set $L_\hi \eqdef +\infty$.

 \section{Proposed Algorithm}

 The proposed \algno is shown as Algorithm~\ref{algo1}. It solves the problem \eqref{eqpbm} as follows. At iteration $t\geq 0$, there are 3 steps: 
\begin{enumerate}
\item Prediction: given the current estimate of the solution $x^{t}$ and the average $\vv^t$ of the dual variables $u_i^t$, the new estimate $\hat{x}^t$ is computed by applying a Gradient Descent (GD) step with respect to $\ff$ and the prox of $\g$.
\item Random activation of the $\hi$: a random subset $\Omega^t\subset [n]$ is chosen and the estimate $y_i^{t+1}$ is computed by applying the prox of  $\hi$, for every $i\in\Omega^t$. The dual variables $u_i$, which are subgradients of the $\hi$, are updated using the differences $y_i^{t+1}-\hat{x}^t$.
\item Update of $x$: the new estimate $x^{t+1}$ is the average of the computed $y_i^{t+1}$. However, if $\Omega^t$ is empty and no $y_i^{t+1}$ is computed, $x^{t+1}$ is randomly chosen as $\hat{x}^t$ with probability $\hat{p}$, or $x^t$ with probability $1-\hat{p}$.
\end{enumerate}

 \noindent\textbf{Arbitrary sampling.} The  subset $\Omega^t$ sampled randomly at every iteration $t\geq 0$ follows an \textbf{arbitrary} proper probability distribution $\mathcal{S}$ over the $2^n$ subsets of $[n]$, characterized by the sets of probabilities $(p_i)_{i=1}^n$ and $(\tilde{p}_i)_{i=1}^n$. 
 For every $i\in[n]$, we define $p_i$ as the probability that $i\in \Omega$ for a random set $\Omega\sim\mathcal{S}$. We assume that  $\mathcal{S}$ is proper; that is, $p_i>0$ for every $i\in[n]$. We define $p_\emptyset \in [0,1)$ as the probability  that  a random set $\Omega\sim\mathcal{S}$ is the empty set $\emptyset$.  
We also define the probabilities $(\tilde{p}_i)_{i=1}^n \in (0,1]^n$ with $\sum_{i=1}^n \tilde{p}_i=1$ such that, for every $(v_i)_{i=1}^n \in \mathbb{R}^n$, we have
\begin{align}
\mathbb{E}_{\Omega\sim\mathcal{S}:\Omega\neq \emptyset}\left[ \frac{1}{|\Omega|}\sum_{i \in \Omega} v_i \right]&=\sum_{i=1}^n \tilde{p}_i v_i,\label{eqprob}
\end{align}
where $|\Omega|$ denotes the size of $\Omega$. These probabilities exist because $\frac{1}{|\Omega|}\sum_{i \in \Omega} v_i$ is a convex combination of the  $v_i$ and a convex combination (when taking the expectation) of a convex combination is a convex combination, and its weights do not depend on the $v_i$.
 Moreover, every $\tilde{p}_i$ is positive because $\mathcal{S}$ is proper. In many cases, the probabilities $\tilde{p}_i$ have a closed form expression. For instance,
 \begin{itemize}
 \item $\Omega\sim\mathcal{S}$ is a subset of size $s$ chosen uniformly at random, for a given minibatch size $s\in [n]$. Then $\tilde{p}_i \equiv \frac{1}{n}$ and $p_i \equiv \frac{s}{n}$.
 \item$\Omega\sim\mathcal{S}$ is of size 1. Then $\tilde{p}_i = p_i$ is the probability that $\Omega=\{i\}$. 
 \end{itemize}
 
It is important to note that  the aggregation process forming $x^{t+1}$ is a simple average, which  is blind to the parameters, notably the properties of the functions $\hi$ and the sampling process $\mathcal{S}$, that characterize how the $y_i^{t+1}$ are computed. So, in the distributed setting considered in Section~\ref{secdi1}, where parallel clients compute the $y_i^{t+1}$, the server averages the received $y_i^{t+1}$, without further consideration. Also, $x^{t+1}$ is a \textbf{biased} estimate of the average of all $n$ $\mathrm{prox}_{\hi}$, in contrast to unbiased randomization in existing algorithms such as 
\algn{Point-SAGA}, 
 \algn{Proxskip}, \algn{RandProx}. 
 This makes our extension to arbitrary sampling far from straightforward.

\setlength{\textfloatsep}{14pt}
 
 \begin{figure}[!t]	
\begin{algorithm}[H]
		\caption{\algno}\label{algo1}
		\begin{algorithmic}[1]
			\STATE  \textbf{parameters:} initial points $x^0,u_1^0,\ldots,u_n^0\in\mathcal{X}$, $\vv^0\coloneqq\frac{1}{n} \sum_{i=1}^n u_i^0$;
			probability $\hat{p}\in[0,1]$; sampling distribution $\mathcal{S}$; 
			stepsizes $(\gamma_t)_{t\geq 0}$, $(\eta_i)_{i=1}^n$
			\FOR{$t=0, 1, \ldots$}
			\STATE $\hat{x}^{t} \coloneqq  \mathrm{prox}_{\gamma_t \g}\big(x^t -\gamma_t \nabla \ff(x^t) - \gamma_t \vv^t\big)$
			\STATE sample $\Omega^t\sim \mathcal{S}$ 
			\IF{$\Omega^t\neq \emptyset$}
			\FOR{$i\in\Omega^t$}
			\STATE $y_{i}^{t+1}\coloneqq \mathrm{prox}_{\gamma_t  \eta_i\hi} (\hat{x}^{t} + \gamma_t  \eta_i u_{i}^t )$
			\STATE $u_{i}^{t+1}\coloneqq u_{i}^t+\frac{1}{\gamma_t \eta_i}  \big(\hat{x}^{t}- y_{i}^{t+1}\big)\in \partial  \hi(y_{i}^{t+1})$ 
			\ENDFOR
			\STATE $u_{i}^{t+1}\coloneqq u_{i}^t $,\ $\forall i\in[n]\backslash \Omega^t$ 
			\STATE  $x^{t+1} \coloneqq \frac{1}{|\Omega^t|}\sum_{i \in \Omega^t}  y_{i}^{t+1}$ 
			\STATE  $\vv^{t+1}\coloneqq 
			\frac{1}{n}\sum_{i=1}^nu_i^{t+1}$
			\ELSE
			\STATE $x^{t+1}\coloneqq \left\{\begin{array}{l} \hat{x}^t\ \mbox{with probability}\ \hat{p}\\
			x^t\ \mbox{with probability}\ 1-\hat{p}
			\end{array}\right.$
			\STATE$u_{i}^{t+1}\coloneqq u_{i}^t$,\ $\forall i\in[n]$ 
			\ENDIF
			\ENDFOR
		\end{algorithmic}
		\end{algorithm}\end{figure}

\subsection{Related Work}\label{secsota}

To minimize a sum of smooth functions, \algn{SGD}-type methods, that make calls to stochastic gradient estimates, are perfectly appropriate and form the de facto standard in ML and AI. Moreover, modern randomized optimization algorithms, developed over the last 10 years or so, are variance-reduced: stabilization or compensation mechanisms using auxiliary variables are employed, so that the random noise decreases and eventually vanishes over the iterations \citep{han19,gor202,gow20a}. 
 Thus, variance-reduced versions of \algn{SGD}, such as \algn{SAGA}  \citep{def14}, \algn{SVRG} \citep{joh13,zha13,xia14,hof15,kov202} and \algn{SARAH} \citep{ngu17,li20b,li21}, converge to the exact solution and are orders of magnitude more efficient than their deterministic counterparts.
 
 However, to deal with nonsmooth functions, an algorithm has to make calls to their proxs instead of their gradients, which may not even be defined. Algorithms that make calls to proxs  are called proximal (splitting) algorithms. They make it possible to solve  a broad class of large-scale optimization problems involving nonsmooth functions  \citep{com10,bot14,par14,kom15,bec17,con23,com21,con22}. But these deterministic algorithms are often too slow, and it is important to design \emph{variance-reduced stochastic proximal algorithms}. The literature on this topic is scarce. Recently, the Stochastic Variance-Reduced Proximal Point Method (\algn{SVRP}), a proximal counterpart of   \algn{SVRG}, has been proposed \citep{kha23,tra24}.
On the other hand,  \algn{Point-SAGA} \citep{def16} is the proximal counterpart of \algn{SAGA}. 
It is recovered as a particular case of \algno when $\ff=\g=0$ and uniform sampling of one function is used, i.e.\ $|\Omega^t| \equiv1$,  $p_i=\tilde{p_i}\equiv \frac{1}{n}$. \algn{Point-SAGA} was extended to an arbitrary fixed minibatch size in \citet{con24}. Even for \algn{Point-SAGA}, our results, notably  importance sampling in Corollary \ref{cor1isa}, linear convergence under similarity in Theorem \ref{theosi1} and accelerated sublinear convergence in Theorem \ref{theoac}, are new. Thus, \algno extends \algn{Point-SAGA}  to the more general problem \eqref{eqpbm} with $\ff$ and $\g$, and to arbitrary sampling with possibly $p_\emptyset >0$.

 \algno is a stochastic version of a parallel \algn{Davis--Yin algorithm}. That is, in the full batch setting where $\Omega^t\equiv[n]$, 
 so that all $\mathrm{prox}_{\hi}$ are called at every iteration, \algno is deterministic and iterates
 \begin{equation*}
 \textstyle
 \left\lfloor\begin{array}{l}
\hat{x}^{t} \coloneqq  \mathrm{prox}_{\gamma_t \g}\big(x^t -\gamma_t \nabla \ff(x^t) -\frac{ \gamma_t }{n}\sum_{i=1}^n u_i^{t+1}\big)\\
			y_{i}^{t+1}\coloneqq \mathrm{prox}_{\gamma_t  \eta_i\hi} (\hat{x}^{t} + \gamma_t  \eta_i u_{i}^t ),\ \forall i\in[n]\\
			 u_{i}^{t+1}\coloneqq u_{i}^t+\frac{1}{\gamma_t \eta_i}  \big(\hat{x}^{t}- y_{i}^{t+1}\big),\ \forall i\in[n]\\
			x^{t+1} \coloneqq \frac{1}{n}\sum_{i =1}^n  y_{i}^{t+1}.
 \end{array}\right.
 \end{equation*}
 We recognize a parallel version of the \algn{Davis--Yin algorithm}, a.k.a.\ three-operator splitting algorithm \citep{dav17}, see \citet[Section 9]{con23} and \citet{con22}. Stochastic versions of this algorithm have been proposed with $\nabla \ff$  replaced by a stochastic estimate \citep{yur16,cev18,ped19,yur21}, which is different from our setting of randomly sampling the proxs. 
 
 In addition to the deterministic  \algn{Davis--Yin algorithm} and  \algn{Point-SAGA}, other existing algorithms are recovered as particular cases of \algno, as summarized in Table~\ref{tab1}.

\begin{table}[t]
\caption{Particular cases of \algno}\label{tab1}
\centering
\begin{tabular}{|l|l|}
Method&Assumptions\\
\hline
Point-SAGA&  $\ff=0$, $\g=0$, $|\Omega^t|=1$, $p_i=\frac{1}{n}$\\
Minibatch Point-SAGA&  $\ff=0$, $\g=0$, $|\Omega^t|=s$, $p_i=\frac{s}{n}$\\
RandProx-Minibatch&$|\Omega^t|=s$, $p_i=\frac{s}{n}$, $\mu_\hi=0$\\
Sto. Decoupling Method& $|\Omega^t|=1$, $\mu_{\g}=0$, $\mu_\hi=0$\\
Proxskip& $n=1$, $\g=0$, $\mu_{\dr{h_1}}=0$\\
RandProx-Skip& $n=1$, $\mu_{\dr{h_1}}=0$\\
Davis--Yin algorithm& $n=1$, $|\Omega^t|=1$\\ 
Parallel Davis--Yin alg.& $|\Omega^t|=n$
\end{tabular}
\end{table}

 \algno with a fixed minibatch size $s\in[n]$ and uniform sampling, so $|\Omega^t|\equiv s$, $p_i\equiv\frac{s}{n}$, $\tilde{p}_i\equiv\frac{1}{n}$, reverts to  \algn{RandProx-Minibatch}, an instance of the general  \algn{RandProx} framework \citep{con23rp} in which the prox of a function is replaced by a stochastic \emph{unbiased} estimate. Linear convergence of  \algn{RandProx-Minibatch} was shown only when  $\mu_\ff>0$ or  $\mu_\g>0$; that is, strong convexity of the $\hi$ is not exploited and \algn{Point-SAGA} cannot be analyzed in this context, for instance. Thus, \algno extends \algn{RandProx-Minibatch} to arbitrary sampling, which departs from the  \algn{RandProx} framework since the update of $x$ becomes biased, and exploits strong convexity of the $\hi$.

\algno with $|\Omega^t|\equiv 1$, so that  exactly one function is activated at every iteration, reverts to the Stochastic Decoupling Method (\algn{SDM}) proposed in \citet{mis19} and analyzed only when $\mu_{\g}=0$ and $\mu_{\hi}\equiv 0$. Our result on \algn{SDM} with importance sampling in Corollary \ref{cor1isa} is new and corrects a mistake in \citet{mis19}.
Thus, \algno extends \algn{SDM} to any minibatch size and more generally to arbitrary sampling, to the important case $p_\emptyset >0$, and exploits strong convexity of $\g$ and the $\hi$.

  \algno with $n=1$ and $\mu_{\dr{h_1}}=0$ becomes \algn{RandProx-Skip} \citep{con23rp}, which reverts itself to \algn{Proxskip} \citep{mis19} if $\g=0$. In the case $n=1$, the important parameter is $p_\emptyset$, and exploiting strong convexity of the $\hi$ when $p_\emptyset>0$ is a crucial part of our new analysis.

\section{Linear Convergence Results}

Under strong convexity of $\ff$ or $\g$ or the $\hi$ and additional assumptions,  \algno  converges linearly. Our first result assumes that the functions $\hi$ are smooth.

\begin{theorem}[Linear convergence with strong convexity and smooth $\hi$]\label{theo1}
Suppose that $L_\hi<+\infty$ for every $i\in[n]$. 
In \algno, suppose that $\gamma_t\equiv \gamma$ for some $0<\gamma< \frac{2}{L_{\ff}}$ (or just $\gamma>0$ if $\ff=0$) and
 let a constant  $\hat{\mu}_\h$ be such that
$
 0\leq \hat{\mu}_\h\leq\min_{i\in [n]}  \frac{2 \eta_i \mu_{\hi} L_\hi}{L_\hi+\mu_{\hi}}
$. 
Suppose that $\mu_{\ff}>0$ or $\mu_{\g}>0$ or $\hat{\mu}_\h>0$. Then the solution $x^\star$ of \eqref{eqpbm} exists and is unique. 
Also, suppose that $\hat{p}\leq \frac{1}{1+\gamma \hat{\mu}_\h}$, let 
$
\bar{p}\eqdef p_\emptyset\hat{p}(1+\gamma \hat{\mu}_\h),
$
and suppose that for every $i\in[n]$,
$
\eta_i= \frac{1- p_\emptyset + \bar{p}}{n \tilde{p}_i (1-p_\emptyset)}
$. 
Define the Lyapunov function, for every $t\geq 0$,
\begin{align}
\Psi^t&\eqdef(1+\gamma \hat{\mu}_\h)\sqn{x^{t}-x^\star} +\frac{1- p_\emptyset + \bar{p}}{n}\sum_{i=1}^n \frac{1}{p_i}
\left(\gamma^2  \eta_i + 
\frac{2\gamma  }{L_{\hi}+\mu_{\hi} } \right)
\sqn{u_i^{t}- \nabla \hi(x^\star)}.\notag
\end{align}
Then \algno converges linearly:  for every $t\geq 0$, 
$\Exp{\Psi^{t}}\leq \rho^t \Psi^0$, 
where 
\begin{align}
\rho&\eqdef\max\left(p_\emptyset(1\!-\!\hat{p})+
(1\!-\!p_\emptyset + \bar{p})\frac{\max(1-\gamma\mu_{\ff},\gamma L_{\ff}-1)^2}{(1+\gamma\mu_{\g})(1+\gamma \hat{\mu}_\h)},1-\min_{i\in[n]}  
\frac{2p_i}{\gamma\eta_i(L_{\hi}+\mu_{\hi})+2}
\right)<1.
\end{align}
Also, $x^t$ converges  to $x^\star$ and  $u_i^t$ converges to $\nabla \hi(x^\star)$ for every $i\in [n]$, 
almost surely.
\end{theorem}
This general linear convergence result has several parameters. In the Appendix, we detail as corollaries several cases with simpler statements, and we discuss the relationship of our linear rates with existing results for particular cases. Notably, we show how to select the probabilities $p_i$ to do importance sampling and balance different $L_\hi$'s, correcting an error in \citet{mis19}.

There are several intertwined parameters, so it is difficult to give a general form of the iteration complexity of \algno,  i.e., the asymptotic number of iterations to achieve $\Exp{\Psi^t}\leq \epsilon$ for any $\epsilon>0$. For example, if $\hat{\mu}_\h=0$ and $\hat{p}=1$, it is
$
\widetilde{\mathcal{O}}\left(
\frac{1}{\gamma(\mu_\ff + \mu_\g)  }+\max_{i\in[n]} \frac{1+\gamma \eta_i L_\hi}{p_i}
\right)
$,
where $\widetilde{\mathcal{O}}(\cdot)=\mathcal{O}\big(\cdot \log(\frac{\Psi^0}{\epsilon})\big)$.

Our second linear convergence result holds for functions $\hi$ that can be nonsmooth but satisfy the following similarity assumption. 

\begin{assumption}[$\delta$-similarity]\label{ass2}
For every $i\in [n]$, $\hi$ is differentiable on $\mathcal{X}$,
\footnote{Differentiability is assumed for simplicity of the notations solely. The functions can be nonsmooth and we can replace the gradients by subgradients in \eqref{eqsim0}, following the example of the anaysis of \algn{Point-SAGA} under similarity in \citet{sad24},
but dealing with set-valued subdifferentials makes all notations and derivations twice as long, since the subgradients are not unique.} 
and  there exists $\delta\geq 0$ such that, for every $(x_i)_{i=1}^n\in\mathcal{X}^n$ and solution $x^\star$ to \eqref{eqincl},
\begin{align}
&\frac{1}{n}\sum_{i=1}^n \left\|\left(\nabla \hi(x_i)-\frac{1}{n}\sum_{j=1}^n\nabla \dr{h_j}(x_j)\right)-\left(\nabla \hi(x^\star)-\frac{1}{n}\sum_{j=1}^n\nabla \dr{h_j}(x^\star)\right)\right\|^2\leq  \frac{\delta^2}{n}\sum_{i=1}^n \sqn{x_i-x^\star}.\label{eqsim0}
\end{align}
We note that if every $\hi$ is $L_\h$-smooth for some $L_\h<+\infty$, \eqref{eqsim0} holds with $\delta=L_\h$, so we can always assume that $\delta\leq L_\h$.
\end{assumption}

\begin{theorem}[Linear convergence with strong convexity and similar $\hi$]\label{theosi1}
Suppose that Assumption~\ref{ass2} holds and, for simplicity, that  $L_\hi\equiv L_\h$ for some $\delta\leq L_\h\leq +\infty$, $\mu_{\hi}\equiv \mu_\h$ for some $\mu_\h\geq 0$, 
$\g=0$, $p_\emptyset=0$, $\tilde{p}_i\equiv\frac{1}{n}$, $p_i\equiv p_s$ for some $p_s\in (0,1]$. In \algno, suppose that $\gamma_t\equiv \gamma$ for some $0<\gamma< \frac{2}{L_{\ff}}$ (or just $\gamma>0$ if $\ff=0$), and $\eta_i\equiv 1$.
Suppose that $\mu_{\ff}>0$ or $\mu_\h>0$. Then the solution $x^\star$ of \eqref{eqpbm} exists and is unique. 
Let $\hat{\mu}_\h\eqdef  \frac{2 \mu_{\h} L_\h}{L_\h+\mu_{\h}}$ (or $2 \mu_{\h}$  if $L_\h=+\infty$) and let $\hat{\mu}_\ff\geq 0$ such that $\max(1-\gamma\mu_{\ff},\gamma L_{\ff}-1)^2=1-\gamma\hat{\mu}_\ff$. 
Suppose that $u_i^0 = \nabla \hi(z_i^0)$ for some $z_i^0\in\mathcal{X}$, for every $i\in [n]$, 
and define the Lyapunov function, for every $t\geq 0$,
\begin{align*}
\Psi^t\!\eqdef\! &\left(1\!-\!\frac{\gamma\hat{\mu}_\ff}{2}\!+\!\frac{\gamma\hat{\mu}_\h}{2} \right)\sqn{x^{t}-x^\star}\!+\!\frac{\gamma(\hat{\mu}_\ff\!+\!\hat{\mu}_\h)}{2n p_s}\sum_{i=1}^n 
\sqn{z_i^{t}-x^\star}+\left(\gamma^2 \! +\! 
\frac{2\gamma  }{L_{\h}+\mu_{\h} } \right)\frac{1}{n p_s}\sum_{i=1}^n
\sqn{u_i^{t}-u_i^\star},
\end{align*}
for auxiliary variables $z_i^t$ defined in Appendix~\ref{apptheosi1}.
Then \algno converges linearly:  for every $t\geq 0$, 
$\Exp{\Psi^{t}}\leq \rho^t \Psi^0$, 
where 
\begin{align}
\rho&\eqdef\max\left(\frac{2-2\gamma\hat{\mu}_\ff}{2-\gamma\hat{\mu}_\ff+\gamma\hat{\mu}_\h},1-p_s\frac{(\hat{\mu}_\ff+\hat{\mu}_\h)(L_\h+\mu_\h)+4\delta^2}{(\hat{\mu}_\ff+\hat{\mu}_\h)(L_\h+\mu_\h)+4\delta^2+2\delta^2\gamma(L_\h+\mu_\h)}\right)<1
\notag\end{align}
(with the second term equal to $1-p_s\frac{(\hat{\mu}_\ff+\hat{\mu}_\h)}{(\hat{\mu}_\ff+\hat{\mu}_\h)+2\delta^2\gamma}$ if $L_\h=+\infty$). 

As a consequence, the iteration complexity of \algno 
is 
$
\widetilde{\mathcal{O}}\left(
\frac{1}{\gamma(\mu_\ff + \mu_\h)  }+\frac{1}{p_s}+
 \frac{\delta^2\gamma }{p_s (\mu_\ff + \mu_\h+\delta^2/L_\h)}
\right)
$.
\end{theorem}

In the smooth heterogeneous case where $\delta\approx L_\h$ 
and there is no
similarity to exploit, the complexity becomes
$
\widetilde{\mathcal{O}}\left(
\frac{1}{\gamma(\mu_\ff + \mu_\h)  }+\frac{1}{p_s}+
 \frac{L^2_\h\gamma }{p_s (\mu_\ff +L_\h)}
\right)
$.
If moreover $\mu_\ff \leq \mu_\h$, which is the case for \algn{Point-SAGA} with $\ff=0$ for instance, the complexity further becomes
$
\widetilde{\mathcal{O}}\left(
\frac{1}{\gamma(\mu_\ff + \mu_\h)  }+\frac{1+L_\h\gamma}{p_s}
\right)
$,
as established for \algn{Point-SAGA}, with $\mu_\ff=0$, in \citet{con24}.

\begin{corollary}[Choosing $\gamma$]
In the conditions of Theorem~\ref{theosi1}, by choosing 
$
\gamma=\min\left(\frac{1}{L_\ff},\frac{\sqrt{p_s}}{\min\big(\delta,\sqrt{L_\h(\mu_\ff + \mu_\h)}\big)}\right)
$
(or by taking a very large $\gamma$ if $\delta=0$ and $L_\ff=0$), 
the complexity of \algno is
\begin{equation}
\widetilde{\mathcal{O}}\left(\frac{L_\ff}{\mu_\ff + \mu_\h}+
\frac{1}{\sqrt{p_s}}\min\!\left(\frac{\delta}{\mu_\ff + \mu_\h},\sqrt{\frac{L_\h}{\mu_\ff + \mu_\h}}
\right)+\frac{1}{p_s}\right).
\end{equation}
\end{corollary}

Thus, \algno converges linearly even without smoothness ($L_\h=+\infty)$, but with smooth $\hi$ and the appropriate value of $\gamma$ it converges even faster with an accelerated complexity depending on 
$\sqrt{\frac{L_\h}{\mu_\ff + \mu_\h}}$. So, \algno converges linearly when the functions $\hi$ are smooth \emph{or} similar, with the best complexity automatically.

In the smooth heterogeneous case where $\delta\approx L_\h$ and $\mu_\ff\leq \mu_\h$, we recover the same result as in Corollary~\ref{cor1usf}. On the other hand, in the nonsmooth case $L_\h=+\infty$, for the particular case of \algn{Point-SAGA} corresponding to $\ff=0$, our new result with a dependence on $\frac{\delta}{\mu_\h}$ is better than the ones in \citet{kha23} and \citet{sad24}, with a dependence on $\frac{\delta^2}{\mu^2_\h}$.

Finally, when $n=1$ and $\g$ is a scaled squared norm, \algno converges linearly \emph{without any assumption on $\dr{h_1}$}. This is consistent with the fact that Assumption~\ref{ass2} is satisfied with $\delta=0$ in this case.

\begin{theorem}[Linear convergence when $n=1$ and $\g$ is simple]\label{theog0}
Suppose that $n=1$, so that $\tilde{p}_1=1$ and $p_1=1-p_\emptyset$, and that $\g=\frac{\mu_\g}{2}\sqn{\cdot}$, for some $\mu_\g\geq 0$.
In \algno, suppose that $\gamma_t\equiv \gamma$ for some $0<\gamma< \frac{2}{L_{\ff}}$ (or just $\gamma>0$ if $\ff=0$) and
 let a constant  $\hat{\mu}_\h$ be such that
$
 0\leq \hat{\mu}_\h\leq \frac{2 \eta_1 \mu_{\dr{h_1}} L_{\dr{h_1}}}{L_{\dr{h_1}}+\mu_{\dr{h_1}}}
$
(or $0\leq \hat{\mu}_\h\leq 2 \eta_1 \mu_{\dr{h_1}}$ if $L_{\dr{h_1}}=+\infty$).
Suppose that $\mu_{\ff}>0$ or $\mu_{\g}>0$ or $\hat{\mu}_\h>0$. Then the solution $x^\star$ of \eqref{eqpbm} exists and is unique. 
Also, suppose that $\hat{p}\leq \frac{1}{1+\gamma \hat{\mu}_\h}$, let 
$
\bar{p}\eqdef p_\emptyset\hat{p}(1+\gamma \hat{\mu}_\h),
$
and suppose that  $\eta_1= \frac{1- p_\emptyset + \bar{p}}{1-p_\emptyset}$. 
Define the Lyapunov function, for every $t\geq 0$,
\begin{align}
\Psi^t\eqdef&(1+\gamma \hat{\mu}_\h)\sqn{x^{t}-x^\star}+
\left(\gamma^2  \eta_1^2 + 
\frac{2\gamma  \eta_1}{L_{\dr{h_1}}+\mu_{\dr{h_1}} } \right)
\sqn{u_1^{t}- \nabla {\dr{h_1}}(x^\star)}.
\end{align}
Then \algno converges linearly:  for every $t\geq 0$, 
$\Exp{\Psi^{t}}\leq \rho^t \Psi^0$, 
where 
{\scriptsize\begin{align*}
\rho&\eqdef\max\!\left(p_\emptyset(1\!-\!\hat{p})+
(1\!-\!p_\emptyset\!+\!\bar{p})\frac{\max(1-\gamma\mu_{\ff},\gamma L_{\ff}-1)^2}{(1+\gamma\mu_{\g})(1+\gamma \hat{\mu}_\h)},
1-\frac{(1-p_\emptyset)^2\big(2(1+\gamma\mu_{\g})+\gamma(L_{\dr{h_1}}+\mu_{\dr{h_1}})\big)
}{\big((1- p_\emptyset + \bar{p})\gamma(L_{\dr{h_1}}+\mu_{\dr{h_1}})+2(1-p_\emptyset)\big)(1+\gamma\mu_{\g})}
\right),
\end{align*}}%
where the second term is equal to 
$
1-\frac{(1-p_\emptyset)^2
}{(1- p_\emptyset + \bar{p})(1+\gamma\mu_{\g})}
$
if $L_{\dr{h_1}}=+\infty$.
Also, $x^t$ converges  to $x^\star$ and  $u_1^t$ converges to $\nabla {\dr{h_1}}(x^\star)$, 
almost surely.
\end{theorem}

\begin{corollary}[$\mu_\hi\equiv 0$]\label{cor7}
In the conditions of Theorem \ref{theog0}, suppose in addition that $\mu_{\dr{h_1}}=\cdots=\mu_{\dr{h_n}}=0$,  so that $\hat{\mu}_\h=0$, 
and that $\hat{p}=1$, so that $\eta_1=\frac{1}{1-p_\emptyset}$. Then
\begin{align*}
\rho=\max&\left(
\frac{\max(1-\gamma\mu_{\ff},\gamma L_{\ff}-1)^2}{1+\gamma\mu_{\g}},
1-\frac{(1-p_\emptyset)^2\big(2(1+\gamma\mu_{\g})+\gamma L_{\dr{h_1}}\big)
}{\big(\gamma L_{\dr{h_1}}+2(1-p_\emptyset)\big)(1+\gamma\mu_{\g})}
\right),
\end{align*}
where the second term is equal to 
$
1-\frac{(1-p_\emptyset)^2
}{1+\gamma\mu_{\g}}
$
if $L_{\dr{h_1}}=+\infty$.
\end{corollary}
In the conditions of Corollary \ref{cor7}, if $\g=0$, so that $\mu_\g=0$, \algno reverts to \algn{ProxSkip} presented in \citet{mis22}.  \algn{ProxSkip} is 
an instance of \algn{RandProx-Skip} \citep{con23rp} and we recover the exact same rate as in \citet[Theorem 2]{con23rp}.  More generally, we see that in the setting  of Corollary \ref{cor7}, it is better to get strong convexity from $\ff$ than from $\g$.

\section{Sublinear Convergence Results}

A $\mathcal{O}(1/t)$ convergence result in the general convex case is shown in the Appendix.

Under strong convexity of $\ff$ or $\g$ or the $\hi$ with no other assumption on the functions, we show accelerated $\mathcal{O}(1/t^2)$ convergence of \algno.
 For simplicity, we distinguish two cases $p_\emptyset=0$ and $\mu_\hi\equiv 0$ (we can deal with the general case for instance by letting $\hat{p}_t$ vary as $\hat{p}_t = \frac{1}{1+\gamma_t \hat{\mu}_\h}$).

\begin{theorem}[$\mathcal{O}(1/t^2)$ convergence, case $p_\emptyset=0$]\label{theoac}
Suppose that in \algno, $p_\emptyset =0$ and $\eta_i=\frac{1}{n\tilde{p}_i}$ for every $i\in[n]$, and that 
$\mu_{\ff}>0$ or $\mu_{\g}>0$ or $\hat{\mu}_\h\eqdef\min_{i\in [n]}  2 \eta_i \mu_{\hi}>0$. Then the solution $x^\star$ of \eqref{eqpbm} exists and is unique.  Let $u_i^\star \in \partial\hi(x^\star)$, $i=1,\ldots,n$, satisfy \eqref{eqincl}. 
Let 
$\mu\eqdef  \max \left(\mu_\ff,\frac{\mu_\g}{2},\frac{\hat{\mu}_\h}{2}\right)>0$
and suppose that $\gamma_t = \frac{2}{\mu (a+t)}$ for every $t\geq 0$, for some constant $a>5$ such that  $\gamma_0\leq \frac{2}{ L_\ff+\mu_\ff}$. 
Define  for every $t\geq 0$ the Lyapunov function
$
\Psi^t\eqdef(1+\gamma_{t-1} \hat{\mu}_\h)\sqn{x^{t}-x^\star} +\frac{1}{n}\sum_{i=1}^n \frac{\gamma_{t-1}^2  \eta_i }{p_i}
\sqn{u_i^{t}- u_i^\star}
$, 
with $\gamma_{-1} \eqdef \frac{2}{\mu(a-1)}$. 
Then for every $t\geq 0$,
\begin{align*}
\Exp{ \sqn{x^t-x^\star}}&\leq  \Exp{\Psi^{t}}\leq \frac{\gamma_{t-1}^2}{\gamma_{-1}^2}\Psi^0=
\frac{(a-1)^2}{(a+t-1)^2}\Psi^0=\mathcal{O}(1/t^2).
\end{align*}
\end{theorem}
To the best of our knowledge, even for particular cases such as \algn{Point-SAGA}, this accelerated convergence result with arbitrary convex functions $\hi$ is new.

\begin{theorem}[$\mathcal{O}(1/t^2)$ convergence, case $\mu_{\hi}\equiv 0$]\label{theoac3}
Suppose that $\mu_{\hi}\equiv 0$, and that $\mu_{\ff}>0$ or $\mu_{\g}>0$.  Then the solution $x^\star$ of \eqref{eqpbm} exists and is unique.  Let $u_i^\star \in \partial\hi(x^\star)$, $i=1,\ldots,n$, satisfy \eqref{eqincl}. Suppose that in \algno,  $\hat{p}=1$ and $\eta_i=\frac{1}{n\tilde{p}_i(1-p_\emptyset)}$ for every $i\in[n]$. 
Let 
$\mu\eqdef  \max \left(\mu_\ff,\frac{\mu_\g}{2}\right)>0$
and suppose that $\gamma_t = \frac{2}{\mu (a+t)}$ for every $t\geq 0$, for some constant $a>5$ such that  $\gamma_0\leq \frac{2}{ L_\ff+\mu_\ff}$. 
Define  for every $t\geq 0$ the Lyapunov function
$
\Psi^t\eqdef \sqn{x^{t}-x^\star} +\frac{1}{n}\sum_{i=1}^n \frac{\gamma_{t-1}^2  \eta_i }{p_i}
\sqn{u_i^{t}- u_i^\star}
$, 
with $\gamma_{-1} \eqdef \frac{2}{\mu(a-1)}$. 
Then for every $t\geq 0$,
\begin{align*}
\Exp{ \sqn{x^t-x^\star}}&\leq \Exp{\Psi^{t}}\leq\frac{\gamma_{t-1}^2}{\gamma_{-1}^2}\Psi^0=
\frac{(a-1)^2}{(a+t-1)^2}\Psi^0=\mathcal{O}(1/t^2).
\end{align*}
\end{theorem}

\section{Distributed Setting}\label{secdi1}

In this section, we consider solving the problem \eqref{eqpbm} in the distributed client-server setting where $n$ clients compute in parallel and communicate with a server. This is particularly relevant in the modern framework of federated learning (FL), where a global model  is trained in a distributed and collaborative way on a large number of machines, such as hospital workstations or mobile phones, where the private data to be exploited is stored  \citep{kon16a,kon16,mcm17, bon17}. 
In this setting, the functions 
$\ff$ and $\g$ are known by the server, which applies $\nabla \ff$ and $\mathrm{prox}_{\g}$, whereas $\hi$ is the private loss function of client $i$, which applies $\mathrm{prox}_{\hi}$. Importantly, we consider the \emph{heterogeneous} setting in which the functions $\hi$ can be arbitrarily different from each other. 
We assume that $\mathcal{X}=\mathbb{R}^d$ for some dimension $d\geq 1$.

In this context, sampling corresponds to \textbf{partial participation}, a.k.a.\ client sampling, an important feature especially in cross-device  FL. Indeed, at a given time, there are multiple reasons for which a a client might be idle and unable to compute or communicate. Partial participation has been well studied for SGD-type methods \citep{gow19,con22mu,con23tam}. Regarding proximal methods, \algn{RandProx-minibatch}, a particular case of \algno as discussed above, allows for partial participation, but in the restricted setting where the cohort size $|\Omega^t|$ 
 is a fixed number $s\geq 1$  chosen randomly by the server. In practice, however, this is not the server deciding whether or not a client participates, the participation of each client depends on its individual choices and constraints. In this context, it is important to allow for arbitrary sampling, as \algno does. In particular, the set $\Omega^t$ of active clients can be empty at a given time $t$,  an important feature of our framework. For example, if every client $i$ participates randomly with probability $p_i>0$, independently of each other, then $p_\emptyset = \prod_{i=1}^n (1-p_i)$. The \algn{5GCS} algorithm proposed by \citet{mal22b} is a variant of \algn{RandProx-minibatch} with inexact computation of the proxs.

In centralized computing, the CPUs or GPUs working in parallel are connected by high performance wire cables or optical fibers. By contrast, in FL, 
communication happens over the internet or mobile phone network and can be slow, costly and unreliable. Thus, communication forms the main bottleneck in FL \citep{kai19,li20,wan21}. To overcome this barrier, \textbf{compression} is used, so that  messages of short bit length are communicated instead of full $d$-dimensional vectors \citep{bez20,alb20,hor22}. \algn{DIANA}  is a distributed GD algorithm that uses compression with a large class of unbiased compressors and is well suited for convex problems \citep{mis19d,hor22d,gor202,con22mu}. 
 \algn{ADIANA} is based on Nesterov-accelerated GD instead \citep{zli20,yhe23}. 
In nonconvex settings,  \algn{MARINA} \citep{marina21} and \algn{DASHA} \citep{dasha23} have been proposed. Some algorithms allow for  biased compressors, such as \algn{EF21} \citep{ric21,fat21} and \algn{EF-BV} \citep{con22e}, but they are less understood, and to date  the best complexity is obtained with unbiased compression.

Apart from compression, a popular strategy to reduce communication is to decrease its frequency. For instance,  communication happens every $10$ or $100$ iterations, or it is triggered randomly with a small probability. This local training approach is at the heart of the popular \algn{FedAvg}  \citep{mcm17} and \algn{Scaffold} \citep{kar20} algorithms and has been analyzed extensively \citep{kha20a,sti19,woo20,li20a,mal20,cha21,gor21,gla22}. 
 \algn{Scaffnew} \citep{mis22}  was the first algorithm exhibiting acceleration from local training. 
Local training and compression have been combined  \citep{bas20,rei20,had21,con24loco}. 
Notably, \citet{con22cs} proposed \algn{CompressedScaffnew}, which harnesses a specific linear compression technique based on permutations
  to \algn{Scaffnew} and achieves double acceleration with respect to the condition number of the functions and the dimension $d$. \algn{TAMUNA} extends this algorithm to partial participation \citep{con23tam}.

\begin{figure}[!t]	
\begin{algorithm}[H]
		\caption{\algnod}\label{algdi}
		\begin{algorithmic}[1]
			\STATE  \textbf{parameters:} initial points $x^0,u_1^0,\ldots,u_n^0\in\mathcal{X}$,  $\vv^0\coloneqq\frac{1}{n} \sum_{i=1}^n u_i^0$;
			probability $\hat{p}\in[0,1]$; sampling distribution $\mathcal{S}$; 
			stepsizes $(\gamma_t)_{t\geq 0}$, $(\eta_i)_{i=1}^n$
			\FOR{$t=0, 1, \ldots$}
			\STATE $\hat{x}^{t} \coloneqq  \mathrm{prox}_{\gamma_t \g}\big(x^t -\gamma_t \nabla \ff(x^t) - \gamma_t \vv^t\big)$\quad// at server
			\STATE sample $\Omega^t\sim \mathcal{S}$ 
			\FOR{$i\in\Omega^t$ at clients,  in parallel}
			\STATE receive $\hat{x}^t$ from server
			\STATE $y_{i}^{t+1}\coloneqq \mathrm{prox}_{\gamma_t  \eta_i\hi} (\hat{x}^{t} + \gamma_t  \eta_i u_{i}^t )$
			\STATE $c_i^t\coloneqq \mathcal{C}_i^t\big(y_{i}^{t+1}-\hat{x}^t\big)$
			\STATE $u_{i}^{t+1}\coloneqq u_{i}^t-\frac{1}{\gamma_t \eta_i} c_i^t$ 
			\STATE send compressed message $c_i^t$ to server
			\ENDFOR
			\STATE $u_{i}^{t+1}\coloneqq u_{i}^t $,\ $\forall i\in[n]\backslash \Omega^t$\quad// idle clients
			\STATE  $x^{t+1} \coloneqq \mathrm{rescale}(\sum_{i\in\Omega^t} c_i^t,\hat{x}^t,x^t)$\ \ // at server, \eqref{eqres}
			\STATE  $\vv^{t+1}\coloneqq \vv^t -\frac{1}{n \gamma_t } \sum_{i\in\Omega^t} \frac{1}{\eta_i} c_i^t$\quad// at server
			\ENDFOR
		\end{algorithmic}
		\end{algorithm}\end{figure}
We present \algnod, shown as Algorithm~\ref{algdi}, which is \textbf{\algno applied to the distributed setting, with the additional feature of \textbf{compression}}. In this algorithm, every client sends a compressed version of $y_i^{t+1}-\hat{x}^t$ to the server. In this paper, we only consider the case where the compressors $\mathcal{C}_i^t$ are independent unscaled \texttt{rand}-$k$ compressors, for some $k\in [d]$. This popular compressor keeps $k$ out of $d$ elements of its argument vector, chosen uniformly at random,  unchanged, and sets the other ones to zero ($\mathcal{C}_i^t(v)$ is the compressed and then decompressed vector, only the $k$ chosen elements are actually communicated, not the zeros, along with a few $\mathcal{O}(k\log d)$ additional bits to indicate the $k$ chosen indexes).

At Step 14 of the algorithm, the $\textbf{rescale}$ procedure averages the received values coordinatewise. That is, let $\textbf{1}$ be the vector of ones of size $d$, $c\eqdef \sum_{i\in\Omega^t} c_i^t$, $s\eqdef \sum_{i\in\Omega^t} \mathcal{C}_i^t(\textbf{1})$ (if $\Omega^t=\emptyset$, $\sum_{i\in\Omega^t}$ is the zero vector). Then for every $j\in [d]$,
\begin{equation}
x^{t+1}_{(j)}\eqdef \left\{\begin{array}{l}\hat{x}^t_{(j)}+c_{(j)}/s_{(j)},\ \mbox{if}\ s_{(j)}>0,\\
\hat{x}^t_{(j)}\ \mbox{with proba.}\ \hat{p}\ \mbox{or}\  
			x^t_{(j)}\ \mbox{with proba.}\ 1-\hat{p},\ \mbox{if}\ s_{(j)}=0,
\end{array}\right.\label{eqres}
\end{equation} 
where $v_{(j)}$ denotes the $j$-th coordinate of the vector $v$. We stress that this adaptive scaling, that takes into account the number of values received by the server for each coordinate, is a central and novel feature of \algnod. By contrast,  in \algn{5GCS-CC} \citep{gru23}, a version of  \algn{5GCS} with compression that also derives from the   \algn{RandProx} framework, \texttt{rand}-$k$ compressors can be used, but the server naively averages the decompressed vectors containing a lot of zeros, to obtain an unbiased estimate. As mentioned above, in \algnod, $x^{t+1}$ is a \textbf{biased} estimate of the average of all $n$ $\mathrm{prox}_{\hi}$, which departs significantly from existing approaches based on unbiased randomization.

The analysis of  \algnod with independent  \texttt{rand}-$k$ compressors follows from the one of \algno, by noticing that the whole analysis is separable in the $d$ coordinates of the vectors, so we can reason about the coordinates independently of each other.  Both sampling and \texttt{rand}-$k$ compression sparsify the vectors, and from the standpoint of a single coordinate, they have the same effect. So, all results on \algno apply, by replacing $p_i$ by $\frac{k}{d}p_i$, 
$p_\emptyset$ by
$
\check{p}_\emptyset\eqdef p_\emptyset+\Exps{\Omega\sim\mathcal{S}}{(1-\frac{k}{d})^{|\Omega|}}
$, 
and the probabilities $(\tilde{p}_i)_{i=1}^n \in (0,1]^n$ are now such that, for every $(v_i)_{i=1}^n \in \mathbb{R}^n$, we have
\begin{align*}
\mathbb{E}_{\Omega\sim\mathcal{S}:\mathcal{C}(\Omega)\neq \emptyset}\left[ \frac{1}{|\mathcal{C}(\Omega)|}\sum_{i \in \mathcal{C}(\Omega)} v_i \right]&=\sum_{i=1}^n \tilde{p}_i v_i,
\end{align*}
where the ``compressed' random set $\mathcal{C}(\Omega)$ is such that for every $i\in\Omega$, $i\in\mathcal{C}(\Omega)$ with probability $\frac{k}{d}$, independently of each other. Note that in general, there is no simple closed form expression of the probabilities $\tilde{p}_i$,
except in the uniform case where all clients behave the same way and $\tilde{p}_i=\frac{1}{n}$, regardless of the other parameters.

We can then consider the (uplink) \textbf{communication complexity}, measured as the number of reals sent in parallel by the clients to the server. This complexity models the communication time and does not depend on the number of clients communicating. 
Thus, the communication complexity is $(1-p_\emptyset)k$ times the iteration complexity.

For instance, as a consequence of Theorem~\ref{theo1}, we have:

\begin{theorem}[Linear convergence with strongly convex and smooth $\hi$, fixed level $s$ of partial participation, $\ff=0$]\label{theod1}
Suppose that $\ff=0$, $L_\hi\equiv L_\h<+\infty$, $\mu_{\hi}\equiv \mu_\h>0$. Then the solution $x^\star$ of \eqref{eqpbm} exists and is unique. 
In \algnod, suppose that $\gamma_t\equiv \gamma>0$ and,
given a batch size $s\in [n]$,  that  $\Omega^t$ is a subset of size $s$ chosen uniformly at random for every $t\geq 0$. Then  $\check{p}_\emptyset=(1-\frac{k}{d})^s$, $\tilde{p}_i \equiv \frac{1}{n}$, $p_i \equiv \frac{ks}{dn}$.
Let
$
 \hat{\mu}_\h\eqdef  \frac{2  \mu_{\h} L_\h}{(1-\check{p}_\emptyset)(L_\h+\mu_{\h})}
$.
Suppose that $\hat{p}= \frac{1}{1+\gamma \hat{\mu}_\h}$ and 
$\eta_i\equiv \frac{1}{1-\check{p}_\emptyset}$.
Define the Lyapunov function, for every $t\geq 0$,
$
\Psi^t\eqdef(1+\gamma \hat{\mu}_\h)\sqn{x^{t}-x^\star} +\frac{d}{ks}
\sum_{i=1}^n 
\left(\gamma^2  \eta_i + 
\frac{2\gamma  }{L_{\h}+\mu_{\h} } \right)
\sqn{u_i^{t}- \nabla \hi(x^\star)}$. 
Then \algnod converges linearly:  for every $t\geq 0$, 
$\Exp{\Psi^{t}}\leq \rho^t \Psi^0$, 
where 
\begin{align*}
\rho&\eqdef\max\left(\frac{1+\check{p}_\emptyset\gamma \hat{\mu}_\h}{1+\gamma \hat{\mu}_\h},1-
\frac{2ks(1-\check{p}_\emptyset)}{dn\big(\gamma(L_{\h}+\mu_{\h})+2(1-\check{p}_\emptyset)\big)}
\right).
\end{align*}
As a consequence, the iteration complexity 
is 
\begin{align}
&\widetilde{\mathcal{O}}\left(
\frac{1}{\gamma\mu_\h  }+\frac{1}{1-\check{p}_\emptyset}+\frac{dn}{ks}+\frac{dn\gamma  L_\h}{ks(1-\check{p}_\emptyset)}
\right).
\end{align}
With $\gamma=\sqrt{\frac{ks(1-\check{p}_\emptyset)}{dn L_\h\mu_\h}}$, it becomes, using $\frac{1}{1-\check{p}_\emptyset}=\mathcal{O}(\frac{d}{ks}+1)$,
$
\widetilde{\mathcal{O}}\left(
\sqrt{\frac{dn L_\h}{ks\mu_\h}}\left(\sqrt{\frac{d}{ks}}+1\right)+\frac{dn}{ks}
\right)
$. 
Consequently, the communication complexity is
$\widetilde{\mathcal{O}}\left(
\frac{d}{s}\sqrt{\frac{n L_\h}{\mu_\h}}+\sqrt{\frac{kdn L_\h}{s\mu_\h}}+\frac{dn}{s}
\right) 
$,
which, as long as $k\leq \max(\frac{d}{s},1)$, simplifies to
\begin{align}
&\widetilde{\mathcal{O}}\left(
\frac{d}{s}\sqrt{\frac{n L_\h}{\mu_\h}}+\sqrt{\frac{dn L_\h}{s\mu_\h}}+\frac{dn}{s}\label{eqcoo1}
\right).
\end{align}
\end{theorem}
The communication complexity of  \algnod in \eqref{eqcoo1} is the same as the one of  
 \algn{TAMUNA} and \algn{5GCS-CC}.

\section{Experiments}\label{secex}

We present a series of computational experiments designed to validate our theoretical findings. In Section~\ref{secex1}, we demonstrate that the importance sampling strategy can significantly outperform its uniform counterpart in terms of convergence rate. Section~\ref{secex2} evaluates the accelerated variant of \algn{Point-SAGA} with adaptive stepsizes -- analyzed in this paper as a special case of \algno{} -- against the classical, non-accelerated version using a fixed stepsize. Finally, in Section~\ref{secex3}, we compare \algnod{} with baseline methods in the distributed setting. All experiments were conducted on a workstation equipped with a 12-core CPU running at 3.7 GHz.

\begin{figure*}
\centering
\includegraphics[width=\linewidth]{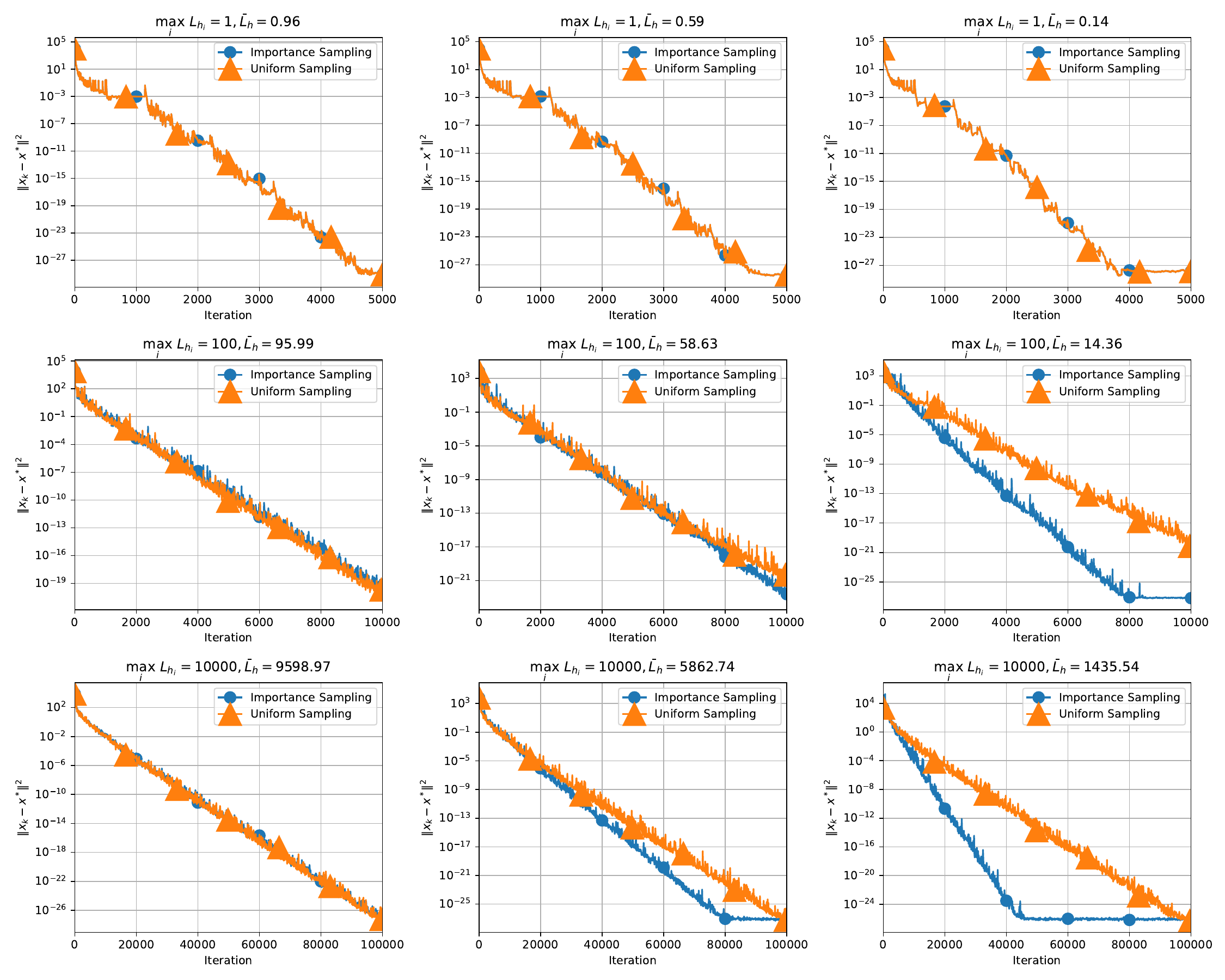}
\caption{Comparison of importance and uniform sampling strategies in \algno. All parameters are set to the theoretical values from Corollaries~\ref{cor1usf} and~\ref{cor1isa}. With uniform sampling, the convergence rate depends crucially on \(\max_i L_{h_i}\), whereas in the importance sampling case, it depends on \(\bar{L}_h = \bigl(\tfrac{1}{n}\sum_{i=1}^n \sqrt{L_{h_i}}\bigr)^2\). When $\bar{L}_h\ll\max_i L_{h_i}$, importance sampling outperforms uniform sampling.}
\label{fig:full_experiment_samplings}
\end{figure*}

\subsection{Importance vs.\ Uniform Sampling}\label{secex1}
Corollary \ref{cor1isa} demonstrates that, in certain scenarios, appropriately tuned importance sampling probabilities yield markedly faster convergence than the uniform sampling baseline established in Corollary \ref{cor1usf}. To test this claim empirically, we consider the following functions:
\[
h_i(x) = \tfrac{1}{2} x^\top \mathbf{A}_i x - \mathbf{b}_i^\top x, 
\quad 
f(x) = \tfrac{1}{2} x^\top \mathbf{B} x, \quad  g \equiv 0,
\]
where each matrix $\mathbf{A}_i \in \mathbb{S}^d_{+}$ is positive semidefinite, $\mathbf{b}_i \in \mathbb{R}^d$ is a vector, and $\mathbf{B} \in \mathbb{S}^d_{+}$ is set to be diagonal for simplicity.  We fix the problem dimensions to $d=100$ and $n=100$, and set the batch size $s = 1$. \algno parameters are selected according to Corollary~\ref{cor1usf} for uniform sampling  and Corollary~\ref{cor1isa} for importance sampling.

To generate each matrix $\mathbf{A}_i$, we first sample a $d$-dimensional orthogonal matrix and a corresponding diagonal matrix, and then combine them according to the spectral decomposition to form $\mathbf{A}_i$.  A fraction of $20\%$ of the diagonal values is set to zero, while the remaining values are drawn uniformly at random from the interval $[0, L_{h_i}]$. In each experiment, 20 functions $h_i$ attain the maximum value $L_{\max} = \max_i L_{h_i}$, while for the remaining 80 functions, $L_{h_i}$ is set to $\alpha L_{\max}$, where $\alpha \in \{0.95, 0.5, 0.05\}$. This variation effectively controls the quantity
$
\bar{L}_h = \left(\frac{1}{n} \sum_{i=1}^{n} \sqrt{L_{h_i}}\right)^2.
$

The diagonal matrix \(\mathbf{B}\) is generated once, with its diagonal entries sampled uniformly at random from the interval \([0.1, 10]\).

Figure \ref{fig:full_experiment_samplings} reveals that the larger the disparity between $\max_i L_{h_i}$ and $\bar{L}_h$, the greater the speedup achieved by importance sampling relative to uniform sampling.

\begin{figure*}
	\centering
	\includegraphics[width=0.7\textwidth]{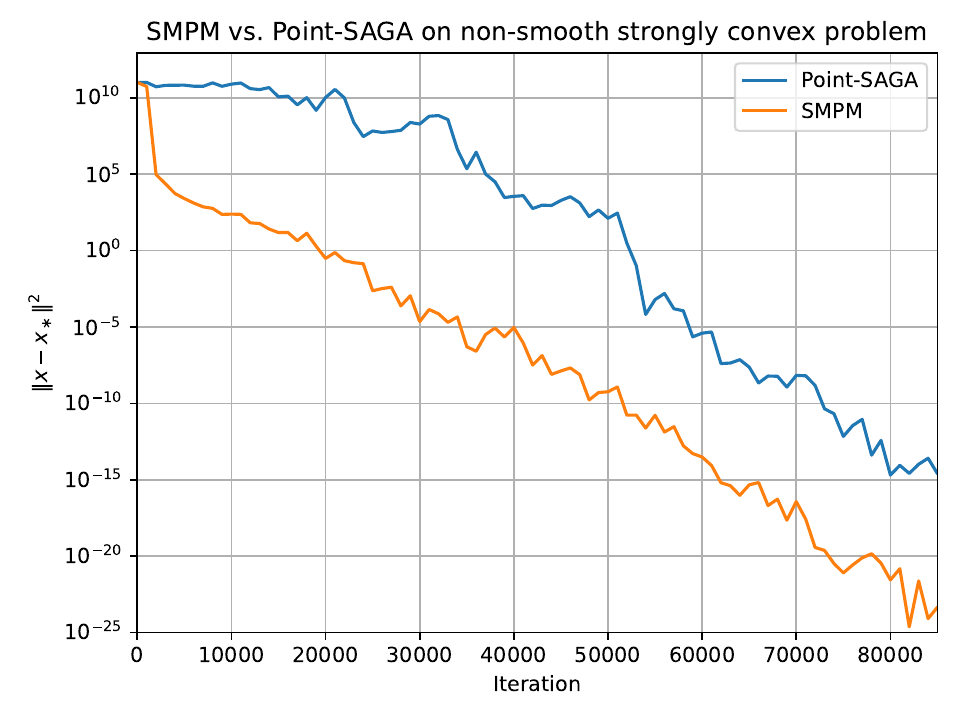}
	\caption{Comparison between \algno and \algn{Point-SAGA}. \algn{Point-SAGA} is run with a fixed stepsize, where the best-performing stepsize is selected from a predefined grid. The configuration of \algno{} follows the setup specified in Theorem~\ref{theoac3}, employing a pre-defined adaptive stepsize schedule that requires no tuning. As shown in the plot, \algno{} consistently outperforms \algn{Point-SAGA}, in line with the theoretical predictions.}
	\label{fig3}
\end{figure*}	

\subsection{\algno vs. \algn{Point-SAGA}}\label{secex2}

As established in Theorem \ref{theoac3}, the adaptive stepsize variant of \algn{Point-SAGA} (see Table \ref{tab1} for the corresponding mapping from \algno{}) attains an accelerated rate compared to the classical constant stepsize schedule in \citep{def16}.  We now verify this theoretical result experimentally.
	
In this experiment, we consider the following setup:
\begin{equation}\label{eq:point_saga_exp_setup}
h_i(x) = \delta(w_i^\top x - b_i) + \frac{\mu}{2} \|x\|^2, \quad f \equiv g \equiv 0, 
\end{equation}
where the indicator function $\delta$ is defined as
$$
\delta(t) = \begin{cases}
0 & \text { if } t = 0, \\
+\infty & \text{ otherwise.}
\end{cases}
$$
Expanding the definition of the proximal operator and applying equation~\eqref{eq:point_saga_exp_setup}, we obtain:
$$\mathrm{prox}_{\gamma h_i}(x)= \frac{x}{\gamma \mu + 1} - \frac{\big(w_i^\top x - b_i (\gamma \mu + 1)\big) w_i}{(\gamma \mu + 1) \|w_i\|^2}.$$
	
We consider $n = 1000$ functions, each defined over a $d = 1000$-dimensional space. The vectors $w_i \in \RR^{d}$ are constructed such that their concatenation forms a randomly generated orthogonal matrix $\bW \in \RR^{n\times d}$. The scalar terms $b_i\in\RR$ are set according to $b = \bW x^\ast$, where $x^\ast\in\RR^d$ is a randomly generated vector. This choice simplifies the computation of the metric $\|x^t - x^\ast\|^2$. The strong convexity parameter $\mu$ is fixed at $10^{-5}$.
	
The stepsize for \algn{Point-SAGA} is set to a constant value  \citep{def16}.
It is selected as the best-performing value from the grid $\{10 \cdot (0.5)^i\}$, where $i \in \{0, \ldots, 11\}$, based on performance over 300,000 iterations. The optimal stepsize within this grid is found to be $10 \cdot 0.5^7$.
	
For \algno{} -- or equivalently, for \algn{Point-SAGA} with adaptive stepsizes -- the stepsize schedule is set according to Theorem~\ref{theoac3} as 
$$
\gamma_t = \frac{2}{\mu (5.5 + t)}.
$$

As shown in Figure~\ref{fig3}, the scheduled stepsize rule outperforms its constant stepsize counterpart. Furthermore, \algno achieves this performance without requiring any stepsize fine-tuning.

\subsection{Comparison between \algnod{}, \algn{5GCS-CC}, and \algn{TAMUNA} in the Distributed Setting} \label{secex3}
	
In this experiment, we empirically compare \algno{} in the distributed setting to its state-of-the-art counterparts, namely, \algn{5GS-CC} and~\algn{TAMUNA}.
	
\begin{figure*}
\centering
\includegraphics[width=0.9\textwidth]{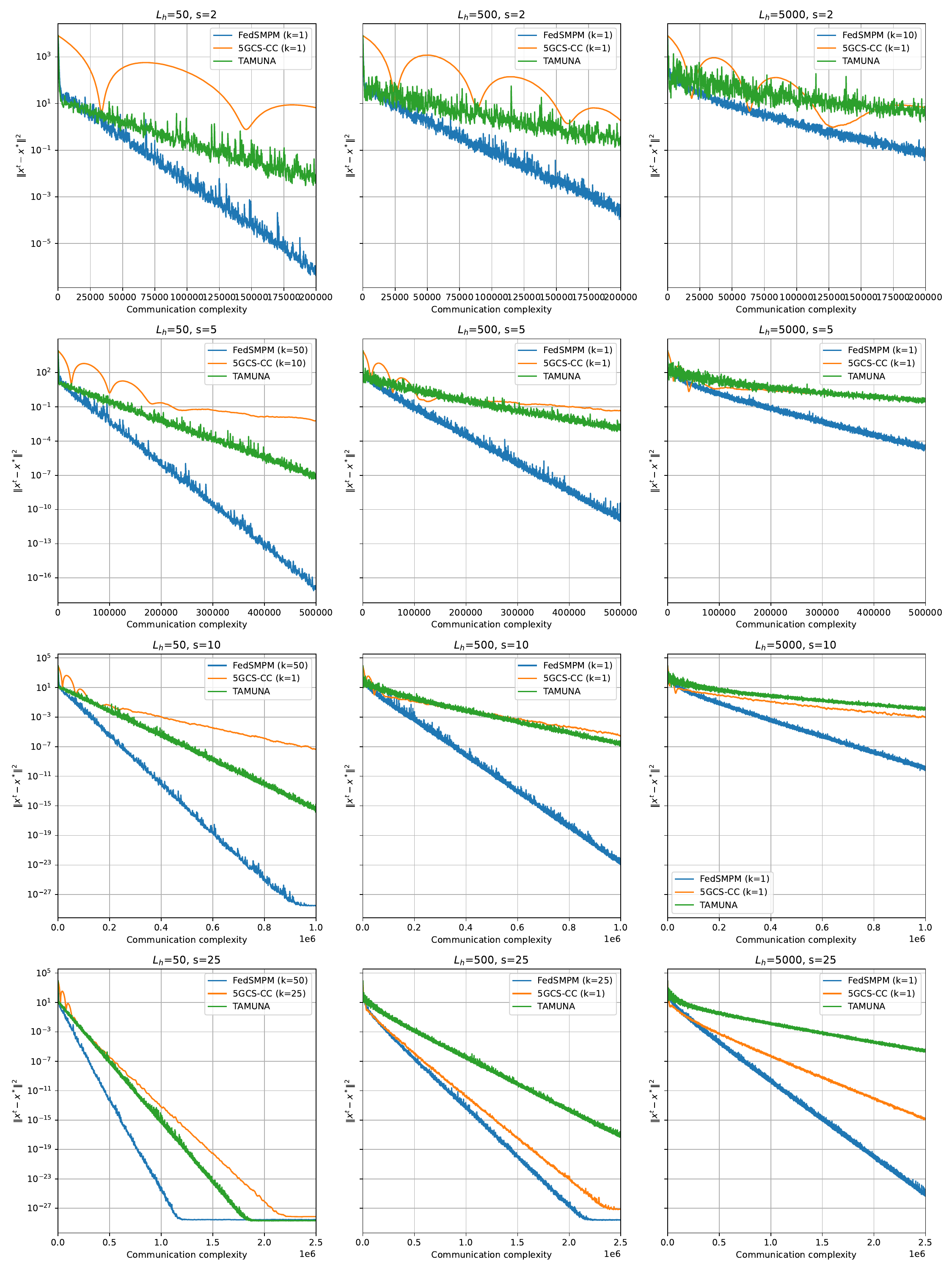}
\caption{Comparison of \algnod{}, \algn{5GCS-CC}, and \algn{TAMUNA} in the distributed setting. \algnod{} and \algn{5GCS-CC} use the \texttt{rand-$k$} compressor, while \algn{TAMUNA} employs the scheme proposed by \citet{con23tam}. All algorithm parameters are set according to their theoretical values. For \algnod{} and \algn{5GCS-CC}, we report the best-performing value of $k$ from the candidate set. In \algn{TAMUNA}, the sparsity parameter $s$ is fixed to $2$. Somewhat unexpectedly, \algnod{} outperforms both baselines.}
\label{figr1}
\end{figure*}
	
We consider the following setup:
$$
h_i(x) = \frac{1}{2} x^\top \mathbf{A}_i x - b_i^\top x, \quad f \equiv g \equiv 0,
$$
where $\mathbf{A}_i \in \mathbb{S}^d$ denotes a symmetric positive definite matrix, and $b_i \in \mathbb{R}^d$ is a vector. In all experiments, we set problem dimensions to $n = 100$ and $d = 100$.
	
The matrix generation process follows the same approach as described in Section~\ref{secex1}. Specifically, to construct each matrix $\mathbf{A}_i$, we first sample their eigenvalues independently from the interval $[\mu, L_{\max}]$, where $\mu$ is fixed to $1$ and $L_{\max} \in \{50, 500, 5000\}$. We then generate independently orthogonal matrices and form each $\mathbf{A}_i$ via  spectral decomposition formula by combining the sampled eigenvalues with the corresponding orthogonal matrix.
To generate $b_i \in \RR^d$, each entry is sampled uniformly from the interval $[0, L_{\max}]$.
	
In all experiments, the starting point is set to  $x_0 = (10, 10, \dots, 10) \in \mathbb{R}^d$.

We implement the \algnod{} algorithm as outlined in Algorithm~\ref{algdi} (see Section~\ref{secdi1}). The algorithm parameters are set according to the theoretical values provided in Theorem~\ref{theod1}.
	
We also replicate the \algn{5GCS-CC} and \algn{TAMUNA} algorithms. The parameters for \algn{5GCS-CC} are selected in accordance with  \citet[Corollary 4.2]{gru23}. For \algn{TAMUNA}, the parameters are configured based on \citet[Theorem 3 and Remark 2]{con23tam}.  
	
We employ the \texttt{rand-k} compression operator when running the \algnod and \algn{5GCS-CC} algorithms. \algnod uses the unscaled version of the operator, while \algn{5GCS-CC} employs the scaled variant. The best performing $k$, selected from the set $\{1, 10, 25, 50\}$ is reported in the plot. For \algn{TAMUNA}, we adopt the compression scheme proposed by \citet{con23tam},  with their sparsity parameter $s$ set to $2$. 
	
As illustrated in Figure~\ref{figr1}, \algnod consistently outperforms the existing baselines.

\section*{Acknowledgments} This work was supported by funding from King Abdullah University of Science and Technology (KAUST): 

\noindent i) KAUST Baseline Research Scheme, 

\noindent ii) Center of Excellence for Generative AI, under award number 5940, 

\noindent iii) Competitive Research Grant (CRG) Program, under award number 6460, 

\noindent iv) SDAIA-KAUST Center of Excellence in Data Science and Artificial Intelligence (SDAIA-KAUST AI).

\bibliographystyle{abbrvnat}
\bibliography{IEEEabrv,biblio2}

\newpage
\appendix

{\huge\noindent\textbf{Appendix}}

\section{General Convergence Analysis}\label{secap1}

We begin the analysis of \algno in the general case.  Let $x^\star$ and $u_i^\star \in \partial\hi(x^\star)$, $i=1,\ldots,n$,  form a solution to \eqref{eqincl}, and let a constant $\hat{\mu}_\h$ be such that
$ 0\leq \hat{\mu}_\h\leq\min_{i\in [n]}  \frac{2 \eta_i \mu_{\hi} L_\hi}{L_\hi+\mu_{\hi}}$ (where the fraction is replaced by $2 \eta_i \mu_{\hi}$  if $L_\hi=+\infty$). Let $t\geq 0$.

For every $i\in[n]$, we define
\begin{equation*}
y_{i}^{t+1}\coloneqq \mathrm{prox}_{\gamma_t \eta_i \hi} (\hat{x}^{t} + \gamma_t \eta_i u_{i}^t )
\end{equation*}
and $\tilde{\nabla} \hi (y_i^{t+1})\in \partial \hi (y_i^{t+1})$ such that
\begin{equation*}
y_{i}^{t+1} +\gamma_t \eta_i \tilde{\nabla} \hi (y_i^{t+1})= \hat{x}^{t} + \gamma_t \eta_i u_{i}^t. 
\end{equation*}
Note that $y_{i}^{t+1}$ and $\tilde{\nabla} \hi (y_i^{t+1})$ are actually computed only for $i \in \Omega^t$.
 
In order to take expectations with respect to the random subset $\Omega^t$, we denote by $\mathcal{F}^t$ the $\sigma$-algebra generated by the collection of random variables $(x^0,u^0_1,\ldots,u^0_n,\ldots, x^t,u^t_1,\ldots,u^t_n)$.

Let $i\in[n]$. We have
\begin{align*}
\Expc{\sqn{u_i^{t+1}-u_i^\star}}
&= p_i \sqn{\tilde{\nabla} \hi(y_i^{t+1})-u_i^\star}
+ (1-p_i)\sqn{u_i^{t} -u_i^\star}.
\end{align*}
Moreover,
\begin{equation*}
y_{i}^{t+1}-x^\star+\gamma_t  \eta_i \big(\tilde{\nabla} \hi(y_i^{t+1})-u_i^\star\big) = \hat{x}^{t} -x^\star + \gamma_t  \eta_i (u_{i}^t-u_i^\star),
\end{equation*}
so that
\begin{align}
&\sqn{y_{i}^{t+1}-x^\star+\gamma_t  \eta_i (\tilde{\nabla} \hi(y_i^{t+1})-u_i^\star) }\notag\\
&=\sqn{y_{i}^{t+1}-x^\star}+\gamma_t^2  \eta_i^2 \sqn{\tilde{\nabla} \hi(y_i^{t+1}) - u_i^\star} + 2\gamma_t  \eta_i \langle
y_{i}^{t+1}-x^\star,\tilde{\nabla} \hi(y_i^{t+1}) - u_i^\star
\rangle\notag\\
&= \sqn{\hat{x}^{t} -x^\star+ \gamma_t  \eta_i (u_{i}^t-u_i^\star)}.\label{eqco1}
\end{align}
According to \citep[Lemma 3.11]{bub15}, 
\begin{align*}
 \langle \tilde{\nabla} \hi(y_i^{t+1})-u_i^\star,y_i^{t+1}-x^\star\rangle &   \geq \frac{\mu_{\hi}L_\hi }{L_\hi+\mu_{\hi}} \sqn{y_i^{t+1}-x^\star} + \frac{1}{L_{\hi}+\mu_{\hi} }  \sqn{\tilde{\nabla} \hi(y_i^{t+1})-u_i^\star}
\end{align*}
(where $ \frac{\mu_{\hi}L_\hi }{L_\hi+\mu_{\hi}} =\mu_\hi$ and $\frac{1}{L_{\hi}+\mu_{\hi} } =0$ if $L_\hi=+\infty$), 
so that
\begin{align*}
\left(1+ \frac{2\gamma_t  \eta_i\mu_{\hi} L_\hi}{L_\hi+\mu_{\hi}} \right)\sqn{y_{i}^{t+1}-x^\star}+\left(\gamma_t^2  \eta_i^2 + 
\frac{2\gamma_t  \eta_i}{L_{\hi}+\mu_{\hi} } \right)
\sqn{\tilde{\nabla} \hi(y_i^{t+1}) - u_i^\star} 
&\leq  \sqn{\hat{x}^{t} -x^\star+ \gamma_t  \eta_i (u_{i}^t-u_i^\star)},
\end{align*}
and, using $\hat{\mu}_\h\leq\min_{i\in [n]}  \frac{2\eta_i \mu_{\hi} L_\hi}{L_\hi+\mu_{\hi}}$,
\begin{align}
&(1+\gamma_t \hat{\mu}_\h)\sqn{y_{i}^{t+1}-x^\star} +\left(\gamma_t^2  \eta_i^2 +
\frac{2 \gamma_t  \eta_i}{L_{\hi}+\mu_{\hi} } \right)
\left(\frac{1}{p_i}\Expc{\sqn{u_i^{t+1}-u_i^\star}}-\frac{1-p_i}{p_i}\sqn{u_i^{t} -u_i^\star }\right)\notag\\
&\leq \sqn{\hat{x}^{t} -x^\star+ \gamma_t  \eta_i (u_{i}^t-u_i^\star)}.\label{eqy7}
\end{align}
If $\Omega^t\neq \emptyset$, $x^{t+1}$ is a convex combination of the $y_{i}^{t+1}$, and by convexity of the squared norm, we have
\begin{align*}
\sqn{x^{t+1}-x^\star}&\leq \frac{1}{|\Omega^t|}\sum_{i \in \Omega^t} \sqn{y_{i}^{t+1}-x^\star}.
\end{align*}
Then, by definition of the $\tilde{p}_i$ in \eqref{eqprob},
\begin{align*}
\Expb{\sqn{x^{t+1}-x^\star}}{\mathcal{F}^t, \Omega^t\neq \emptyset}&\leq 
\sum_{i=1}^n \tilde{p}_i \sqn{y_{i}^{t+1}-x^\star}.
\end{align*}
Also,
\begin{align*}
\Expb{\sqn{x^{t+1}-x^\star}}{\mathcal{F}^t, \Omega^t= \emptyset}&=
\hat{p}\sqn{\hat{x}^t-x^\star}+(1-\hat{p})
\sqn{x^t-x^\star}.
\end{align*}
Thus,
\begin{align*}
\Expc{\sqn{x^{t+1}-x^\star}}&\leq p_\emptyset\hat{p}\sqn{\hat{x}^t-x^\star}+p_\emptyset(1-\hat{p})\sqn{x^t-x^\star}+(1- p_\emptyset)
\sum_{i=1}^n \tilde{p}_i \sqn{y_{i}^{t+1}-x^\star}
\end{align*}
and
\begin{align}
\frac{1}{1- p_\emptyset}\left(\Expc{\sqn{x^{t+1}-x^\star}}- p_\emptyset\hat{p}\sqn{\hat{x}^t-x^\star}-p_\emptyset(1-\hat{p})\sqn{x^t-x^\star}\right)&\leq \sum_{i=1}^n \tilde{p}_i\sqn{y_{i}^{t+1}-x^\star}.\label{eqco2}
\end{align}
Combining this inequality with \eqref{eqy7}, we obtain
\begin{align*}
&\frac{1+\gamma_t \hat{\mu}_\h}{1- p_\emptyset}\left(\Expc{\sqn{x^{t+1}-x^\star}}- p_\emptyset\hat{p}\sqn{\hat{x}^t-x^\star}-p_\emptyset(1-\hat{p})\sqn{x^t-x^\star}\right)\\
&\quad +\sum_{i=1}^n  \tilde{p}_i\left(\gamma_t^2  \eta_i^2 + 
\frac{2\gamma_t  \eta_i}{L_{\hi}+\mu_{\hi} } \right)
\left(\frac{1}{p_i}\Expc{\sqn{u_i^{t+1}-u_i^\star}}-\frac{1-p_i}{p_i}\sqn{u_i^{t} -u_i^\star }\right)\\
&\leq (1+\gamma_t \hat{\mu}_\h)\sum_{i=1}^n \tilde{p}_i\sqn{y_{i}^{t+1}-x^\star}\\
&\quad +\sum_{i=1}^n  \tilde{p}_i\left(\gamma_t^2  \eta_i^2 +
\frac{2\gamma_t  \eta_i}{L_{\hi}+\mu_{\hi} } \right)
\left(\frac{1}{p_i}\Expc{\sqn{u_i^{t+1}-u_i^\star}}-\frac{1-p_i}{p_i}\sqn{u_i^{t} -u_i^\star }\right)\\
&\leq  \sum_{i=1}^n\tilde{p}_i \sqn{\hat{x}^{t} -x^\star+ \gamma_t  \eta_i (u_{i}^t-u_i^\star)},
\end{align*}
so that
\begin{align}
&\frac{1+\gamma_t \hat{\mu}_\h}{1- p_\emptyset}\Expc{\sqn{x^{t+1}-x^\star}} +\sum_{i=1}^n \frac{\tilde{p}_i}{p_i}
\left(\gamma_t^2  \eta_i^2 + 
\frac{2\gamma_t  \eta_i}{L_{\hi}+\mu_{\hi} } \right)
\Expc{\sqn{u_i^{t+1}-u_i^\star}}\notag\\
&\leq \frac{p_\emptyset \hat{p}(1+\gamma_t \hat{\mu}_\h)}{1- p_\emptyset} \sqn{\hat{x}^{t}-x^\star}+\frac{p_\emptyset(1-\hat{p})(1+\gamma_t \hat{\mu}_\h)}{1-p_\emptyset}\sqn{x^t-x^\star}\notag\\
&\quad+
\sum_{i=1}^n 
\frac{ \tilde{p}_i(1-p_i)}{p_i}\left(\gamma_t^2  \eta_i^2 + 
\frac{2\gamma_t  \eta_i}{L_{\hi}+\mu_{\hi} } \right)
\sqn{u_i^{t} -u_i^\star }+ \sum_{i=1}^n \tilde{p}_i\sqn{\hat{x}^{t} -x^\star+ \gamma_t  \eta_i (u_{i}^t-u_i^\star)}.\label{eqpp0}
\end{align}
We suppose that $p_\emptyset=0$, or $\hat{\mu}_\h=0$, or $\gamma_t$ is constant, and that 
 $\hat{p}\leq \frac{1}{1+\gamma_t \hat{\mu}_\h}$. We let 
\begin{equation*}
\bar{p}\eqdef p_\emptyset\hat{p}(1+\gamma_t \hat{\mu}_\h)\leq p_\emptyset
\end{equation*}
(that does not depend on $t$ under the previous assumptions), 
and we suppose that for every $i\in[n]$,
\begin{equation*}
\eta_i= \frac{1- p_\emptyset + \bar{p}}{n \tilde{p}_i (1-p_\emptyset)}.
\end{equation*}
Multiplying \eqref{eqpp0} by $1- p_\emptyset$, we obtain
\begin{align*}
&(1+\gamma_t \hat{\mu}_\h)\Expc{\sqn{x^{t+1}-x^\star}} +(1- p_\emptyset + \bar{p})\frac{1}{n}\sum_{i=1}^n \frac{1}{p_i}
\left(\gamma_t^2  \eta_i + 
\frac{2\gamma_t  }{L_{\hi}+\mu_{\hi} } \right)
\Expc{\sqn{u_i^{t+1}-u_i^\star}}\notag\\
&\leq p_\emptyset(1-\hat{p})(1+\gamma_t \hat{\mu}_\h)\sqn{x^t-x^\star}+
(1- p_\emptyset + \bar{p})\frac{1}{n}\sum_{i=1}^n 
\frac{ 1-p_i}{p_i}\left(\gamma_t^2  \eta_i + 
\frac{2\gamma_t  }{L_{\hi}+\mu_{\hi} } \right)
\sqn{u_i^{t} -u_i^\star }\\
&\quad + \bar{p} \sqn{\hat{x}^{t}-x^\star}+ (1- p_\emptyset)\sum_{i=1}^n \tilde{p}_i\sqn{\hat{x}^{t} -x^\star+ \gamma_t  \eta_i (u_{i}^t-u_i^\star)}\\
&= p_\emptyset(1-\hat{p})(1+\gamma_t \hat{\mu}_\h)\sqn{x^t-x^\star}+
(1- p_\emptyset + \bar{p})\frac{1}{n}\sum_{i=1}^n 
\frac{ 1-p_i}{p_i}\left(\gamma_t^2  \eta_i + 
\frac{2\gamma_t  }{L_{\hi}+\mu_{\hi} } \right)
\sqn{u_i^{t} -u_i^\star }\\
&\quad + (1- p_\emptyset + \bar{p})\left(\frac{\bar{p}}{1- p_\emptyset + \bar{p}} \sqn{\hat{x}^{t}-x^\star}+ \frac{1- p_\emptyset}{1- p_\emptyset + \bar{p}}\sum_{i=1}^n \tilde{p}_i\sqn{\hat{x}^{t} -x^\star+ \gamma_t  \eta_i (u_{i}^t-u_i^\star)}\right).
\end{align*}
To develop the last term, let $\tilde{\nabla}\g(\hat{x}^{t})\in\partial \g(\hat{x}^{t})$ be such that 
\begin{equation*}
\hat{x}^{t}+\gamma_t \tilde{\nabla}\g(\hat{x}^{t}) =x^t-\gamma_t \nabla \ff(x^t)-\frac{\gamma_t}{n} \sum_{i=1}^n u_i^t.
\end{equation*}
Let $q^t \coloneqq  \tilde{\nabla}\g(\hat{x}^{t})- \mu_{\g} \hat{x}^{t}$.
Then 
\begin{equation*}
(1+\gamma_t\mu_{\g})\hat{x}^{t}=x^t-\gamma_t \nabla \ff(x^t)-\gamma_t q^t -\frac{\gamma_t}{n} \sum_{i=1}^n u_i^t.
\end{equation*}
Similarly, we define $q^\star$ such that 
$(1+\gamma_t\mu_{\g})x^\star=x^\star-\gamma_t \nabla \ff(x^\star)-\gamma_t q^\star -\frac{\gamma_t}{n} \sum_{i=1}^n u_i^\star$. Also,
 let $w^t\coloneqq x^t-\gamma_t \nabla \ff(x^t)$ and $w^\star\coloneqq x^\star-\gamma_t \nabla \ff(x^\star)$. 
 We have
\begin{align}
&\frac{\bar{p}}{1- p_\emptyset + \bar{p}} \sqn{\hat{x}^{t}-x^\star}+ \frac{1- p_\emptyset}{1- p_\emptyset + \bar{p}}\sum_{i=1}^n \tilde{p}_i\sqn{\hat{x}^{t} -x^\star+ \gamma_t  \eta_i (u_{i}^t-u_i^\star)}\notag\\
&=\sqn{\hat{x}^{t}-x^\star} 
+ 2\gamma_t  \frac{1- p_\emptyset}{1- p_\emptyset + \bar{p}}\left\langle \hat{x}^{t} -x^\star, \sum_{i=1}^n \tilde{p}_i  \eta_i (u_{i}^t-u_i^\star)\right\rangle\notag\\
&\quad+ \gamma_t^2  \frac{1- p_\emptyset}{1- p_\emptyset + \bar{p}}\sum_{i=1}^n \tilde{p}_i  \eta_i^2 \sqn{u_{i}^t-u_i^\star}\notag\\
&=(1+\gamma_t\mu_{\g})\sqn{\hat{x}^{t}-x^\star} -\gamma_t\mu_{\g}\sqn{\hat{x}^{t}-x^\star}\notag\\
&\quad+ 2\gamma_t \left\langle \hat{x}^{t} -x^\star,\frac{1}{n} \sum_{i=1}^n  (u_{i}^t-u_i^\star)\right\rangle
+\frac{ \gamma_t^2}{n} \sum_{i=1}^n   \eta_i \sqn{u_{i}^t-u_i^\star}\notag\\
&= \left\langle w^t-w^\star-\gamma_t (q^t-q^\star) -\frac{\gamma_t}{n} \sum_{j=1}^n (u_j^t-u_j^\star),\hat{x}^{t}-x^\star\right\rangle -\gamma_t\mu_{\g}\sqn{\hat{x}^{t}-x^\star}\notag\\
&\quad+ 2\gamma_t \left\langle \hat{x}^{t} -x^\star,\frac{1}{n} \sum_{i=1}^n  (u_{i}^t-u_i^\star)\right\rangle
+\frac{ \gamma_t^2}{n} \sum_{i=1}^n   \eta_i \sqn{u_{i}^t-u_i^\star} \notag\\
&= \left\langle w^t-w^\star+\gamma_t (q^t-q^\star) +\frac{\gamma_t}{n} \sum_{j=1}^n (u_j^t-u_j^\star),\hat{x}^{t}-x^\star\right\rangle\notag\\
&\quad-\gamma_t\mu_{\g}\sqn{\hat{x}^{t}-x^\star} -2\gamma_t \langle q^t-q^\star, \hat{x}^{t} -x^\star\rangle+\frac{ \gamma_t^2}{n} \sum_{i=1}^n   \eta_i \sqn{u_{i}^t-u_i^\star}\notag\\
&= \frac{1}{1+\gamma_t\mu_{\g}}\left\langle w^t-w^\star+\gamma_t (q^t-q^\star) +\frac{\gamma_t}{n} \sum_{j=1}^n (u_j^t-u_j^\star),\right.\notag\\
&\quad\left. w^t-w^\star-\gamma_t (q^t-q^\star) -\frac{\gamma_t}{n} \sum_{j=1}^n (u_j^t-u_j^\star)\right\rangle -\gamma_t\mu_{\g}\sqn{\hat{x}^{t}-x^\star}\notag\\
&\quad- 2\gamma_t \langle q^t-q^\star, \hat{x}^{t} -x^\star\rangle +\frac{ \gamma_t^2}{n} \sum_{i=1}^n   \eta_i \sqn{u_{i}^t-u_i^\star}\notag\\
&= \frac{1}{1+\gamma_t\mu_{\g}}\sqn{w^t-w^\star}- \frac{\gamma_t^2}{1+\gamma_t\mu_{\g}}\sqn{q^t-q^\star+\frac{1}{n} \sum_{j=1}^n (u_j^t-u_j^\star)}
\notag\\
&\quad- 2\gamma_t \langle q^t-q^\star, \hat{x}^{t} -x^\star\rangle  -\gamma_t\mu_{\g}\sqn{\hat{x}^{t}-x^\star}+\frac{ \gamma_t^2}{n} \sum_{i=1}^n   \eta_i \sqn{u_{i}^t-u_i^\star}.\label{eqco4}
 \end{align}
 Hence,
\begin{align*}
&(1+\gamma_t \hat{\mu}_\h)\Expc{\sqn{x^{t+1}-x^\star}} +\frac{1- p_\emptyset + \bar{p}}{n}\sum_{i=1}^n \frac{1}{p_i}
\left(\gamma_t^2  \eta_i + 
\frac{2\gamma_t  }{L_{\hi}+\mu_{\hi} } \right)
\Expc{\sqn{u_i^{t+1}-u_i^\star}}\notag\\
&\leq  p_\emptyset(1-\hat{p})(1+\gamma_t \hat{\mu}_\h)\sqn{x^t-x^\star}+
\frac{1- p_\emptyset + \bar{p}}{1+\gamma_t\mu_{\g}}\sqn{w^t-w^\star}\\
&\quad+\frac{1- p_\emptyset + \bar{p}}{n}\sum_{i=1}^n \left(
\frac{ 1-p_i}{p_i}\left(\gamma_t^2  \eta_i + 
\frac{2\gamma_t  }{L_{\hi}+\mu_{\hi} } \right)+\gamma_t^2  \eta_i \right)\sqn{u_i^{t} -u_i^\star }\\
&\quad- \frac{\gamma_t^2(1- p_\emptyset + \bar{p})}{1+\gamma_t\mu_{\g}}\sqn{q^t-q^\star+\frac{1}{n} \sum_{j=1}^n (u_j^t-u_j^\star)}
- 2\gamma_t(1- p_\emptyset + \bar{p}) \langle q^t-q^\star, \hat{x}^{t} -x^\star\rangle\\
&\quad-\gamma_t\mu_{\g}(1- p_\emptyset + \bar{p})\sqn{\hat{x}^{t}-x^\star}.
\end{align*}
According to \citep[Lemma 1]{con23rp},
\begin{eqnarray}
\sqn{w^t-w^\star}&=& \sqn{(\mathrm{Id}-\gamma_t\nabla \ff)x^t-(\mathrm{Id}-\gamma_t\nabla \ff)x^\star}\notag \\
&\leq &\max(1-\gamma_t\mu_{\ff},\gamma_t L_{\ff}-1)^2 \sqn{x^t-x^\star}.\label{eqw1}
\end{eqnarray}
Thus,
\begin{align}
&(1+\gamma_t \hat{\mu}_\h)\Expc{\sqn{x^{t+1}-x^\star}} +\frac{1- p_\emptyset + \bar{p}}{n}\sum_{i=1}^n \frac{1}{p_i}
\left(\gamma_t^2  \eta_i + 
\frac{2\gamma_t  }{L_{\hi}+\mu_{\hi} } \right)
\Expc{\sqn{u_i^{t+1}-u_i^\star}}\notag\\
&\leq  \left(p_\emptyset(1-\hat{p})+
(1- p_\emptyset + \bar{p})\frac{\max(1-\gamma_t\mu_{\ff},\gamma_t L_{\ff}-1)^2}{(1+\gamma_t\mu_{\g})(1+\gamma_t \hat{\mu}_\h)}\right)(1+\gamma_t \hat{\mu}_\h)\sqn{x^t-x^\star}\notag\\
&\quad+\frac{1- p_\emptyset + \bar{p}}{n}\sum_{i=1}^n \left(1-p_i + p_i\frac{\gamma_t\eta_i(L_{\hi}+\mu_{\hi})}{\gamma_t\eta_i(L_{\hi}+\mu_{\hi})+2}
\right)
\frac{1}{p_i}
\left(\gamma_t^2  \eta_i + 
\frac{2\gamma_t  }{L_{\hi}+\mu_{\hi} } \right)\sqn{u_i^{t} -u_i^\star }\notag\\
&\quad- \frac{\gamma_t^2(1- p_\emptyset + \bar{p})}{1+\gamma_t\mu_{\g}}\sqn{q^t-q^\star+\frac{1}{n} \sum_{j=1}^n (u_j^t-u_j^\star)}
- 2\gamma_t(1- p_\emptyset + \bar{p}) \langle q^t-q^\star, \hat{x}^{t} -x^\star\rangle\notag\\
&\quad-\gamma_t\mu_{\g}(1- p_\emptyset + \bar{p})\sqn{\hat{x}^{t}-x^\star}\notag\\
&=  \left(p_\emptyset(1-\hat{p})+
(1- p_\emptyset + \bar{p})\frac{\max(1-\gamma_t\mu_{\ff},\gamma_t L_{\ff}-1)^2}{(1+\gamma_t\mu_{\g})(1+\gamma_t \hat{\mu}_\h)}\right)(1+\gamma_t \hat{\mu}_\h)\sqn{x^t-x^\star}\label{eqp1}\\
&\quad+\frac{1- p_\emptyset + \bar{p}}{n}\sum_{i=1}^n \left(1-\frac{2p_i}{\gamma_t\eta_i(L_{\hi}+\mu_{\hi})+2}\right)
\frac{1}{p_i}
\left(\gamma_t^2  \eta_i + 
\frac{2\gamma_t  }{L_{\hi}+\mu_{\hi} } \right)\sqn{u_i^{t} -u_i^\star }\notag\\
&\quad- \frac{\gamma_t^2(1- p_\emptyset + \bar{p})}{1+\gamma_t\mu_{\g}}\sqn{q^t-q^\star+\frac{1}{n} \sum_{j=1}^n (u_j^t-u_j^\star)}
- 2\gamma_t(1- p_\emptyset + \bar{p}) \langle q^t-q^\star, \hat{x}^{t} -x^\star\rangle\notag\\
&\quad-\gamma_t\mu_{\g}(1- p_\emptyset + \bar{p})\sqn{\hat{x}^{t}-x^\star}.\notag
\end{align}

\section{Proof of Theorem~\ref{theo1}}

Let $t\geq 0$.  Following the derivations in Section \ref{secap1}, 
we have from \eqref{eqp1} that
\begin{equation}
\Expc{\Psi^{t+1}}\leq \rho \Psi^t,\label{eqrec2bz}
\end{equation}
where
\begin{align*}
\rho&\eqdef\max\left(p_\emptyset(1-\hat{p})+
(1- p_\emptyset + \bar{p})\frac{\max(1-\gamma\mu_{\ff},\gamma L_{\ff}-1)^2}{(1+\gamma\mu_{\g})(1+\gamma \hat{\mu}_\h)},
1-\min_{i\in[n]}  
\frac{2p_i}{\gamma\eta_i(L_{\hi}+\mu_{\hi})+2}
\right).
\end{align*}
Using the tower rule, we can unroll the recursion in \eqref{eqrec2bz} to obtain the unconditional expectation of $\Psi^{t+1}$. 

Under the assumptions of  Theorem \ref{theo1}, $\rho<1$. Indeed the first term in $\rho$ is
\begin{align*}
&p_\emptyset-p_\emptyset\hat{p}+
(1- p_\emptyset)\frac{\max(1-\gamma\mu_{\ff},\gamma L_{\ff}-1)^2}{(1+\gamma\mu_{\g})(1+\gamma \hat{\mu}_\h)}
+p_\emptyset\hat{p}\frac{\max(1-\gamma\mu_{\ff},\gamma L_{\ff}-1)^2}{(1+\gamma\mu_{\g})}\\
&\leq p_\emptyset +(1- p_\emptyset)\frac{\max(1-\gamma\mu_{\ff},\gamma L_{\ff}-1)^2}{(1+\gamma\mu_{\g})(1+\gamma \hat{\mu}_\h)}<1.
\end{align*}
In words, there are 3 cases: 1) with probability $1-p_\emptyset>0$ (steps 8--14) there is a contraction by $\frac{\max(1-\gamma\mu_{\ff},\gamma L_{\ff}-1)^2}{(1+\gamma\mu_{\g})(1+\gamma \hat{\mu}_\h)}<1$; 2) with probability $p_\emptyset \hat{p}$, $x^{t+1}$ is updated with $\hat{x}$ and there is a contraction by $\frac{\max(1-\gamma\mu_{\ff},\gamma L_{\ff}-1)^2}{(1+\gamma\mu_{\g})}\leq 1$; 3) with probability $p_\emptyset (1-\hat{p})$, $x^{t+1}$ is updated with $x^t$ and there is no contraction.

Since $\rho<1$, $\Exp{\Psi^t}\rightarrow 0$. 
Moreover, using classical results on supermartingale convergence \citep[Proposition A.4.5]{ber15}, it follows from \eqref{eqrec2bz} that $\Psi^t\rightarrow 0$ almost surely. Almost sure convergence of $x^t$ and the $u_i^t$ follows.

\section{Corollaries of Theorem~\ref{theo1}}

First, we consider the case  $p_\emptyset =0$ in which at least one function $\hi$ is activated at every iteration. $\hat{p}$ does not need to be defined in that case, so that $\hat{\mu}_\h$ does not need to be known and we can set $\bar{p}=p_\emptyset=0$.

\begin{corollary}[Case $p_\emptyset =0$]\label{cor1}
Suppose that $L_\hi<+\infty$ for every $i\in[n]$. 
In \algno, suppose that $\gamma_t\equiv \gamma$ for some $0<\gamma< \frac{2}{L_{\ff}}$ (or just $\gamma>0$ if $\ff=0$) and $\eta_i= \frac{1}{n \tilde{p}_i }$ for every $i\in[n]$.
Let \begin{equation}
\hat{\mu}_\h\eqdef \min_{i\in [n]}  \frac{2  \mu_{\hi} L_\hi}{n \tilde{p}_i(L_\hi+\mu_{\hi})}.
\end{equation}
Suppose that $\mu_{\ff}>0$ or $\mu_{\g}>0$ or $\hat{\mu}_\h>0$. Then the solution $x^\star$ of \eqref{eqpbm} exists and is unique. 
Define the Lyapunov function, for every $t\geq 0$,
\begin{equation}
\Psi^t\eqdef(1+\gamma \hat{\mu}_\h)\sqn{x^{t}-x^\star} +\frac{1}{n}\sum_{i=1}^n \frac{1}{p_i}
\left(\gamma^2  \eta_i + 
\frac{2\gamma  }{L_{\hi}+\mu_{\hi} } \right)
\sqn{u_i^{t}- \nabla \hi(x^\star)}.
\end{equation}
Then \algno converges linearly:  for every $t\geq 0$, 
$\Exp{\Psi^{t}}\leq \rho^t \Psi^0$, 
where 
\begin{align}
\rho&\eqdef\max\left(\frac{\max(1-\gamma\mu_{\ff},\gamma L_{\ff}-1)^2}{(1+\gamma\mu_{\g})(1+\gamma \hat{\mu}_\h)},
1-\min_{i\in[n]}  
\frac{2n \tilde{p}_i p_i}{\gamma(L_{\hi}+\mu_{\hi})+2n \tilde{p}_i}
\right)<1.
\end{align}
Also, $x^t$ converges  to $x^\star$ and  $u_i^t$ converges to $\nabla \hi(x^\star)$ for every $i\in [n]$, 
almost surely.

As a consequence, the iteration complexity of \algno is 
\begin{equation}
\mathcal{O}\left(\left(
\frac{1}{\gamma \left(\mu_\ff + \mu_\g+\min_{i\in [n]}  \frac{\mu_{\hi}}{n\tilde{p}_i}\right) }
+\max_{i\in[n]} \left(\frac{1}{p_i}+\frac{\gamma  L_\hi}{n\tilde{p}_i p_i}\right)
\right)\log\left(\frac{\Psi^0}{\epsilon}\right)\right).
\end{equation}
\end{corollary}

\begin{corollary}[Uniform minibatch sampling]\label{cor1us}
Given a batch size $s\in [n]$,  suppose that in \algno $\mathcal{S}$ is a uniform distribution, so that $\Omega^t$ is a subset of size $s$ chosen uniformly at random for every $t\geq 0$. Then  $p_\emptyset=0$, $\tilde{p}_i \equiv \frac{1}{n}$, $p_i \equiv \frac{s}{n}$, and  in the conditions of Corollary \ref{cor1}, we have $\eta_i\equiv 1$ and
\begin{align}
\rho&=\max\left(\frac{\max(1-\gamma\mu_{\ff},\gamma L_{\ff}-1)^2}{(1+\gamma\mu_{\g})(1+\gamma \hat{\mu}_\h)},
1-\frac{s}{n}\min_{i\in[n]}  
\frac{2 }{\gamma(L_{\hi}+\mu_{\hi})+2}
\right)<1.
\end{align}
As a consequence, the iteration complexity of \algno is 
\begin{equation}
\mathcal{O}\left(\left(
\frac{1}{\gamma(\mu_\ff + \mu_\g+\min_{i\in [n]} \mu_{\hi}) }
+\frac{n}{s}+\frac{n}{s}\max_{i\in[n]}\gamma  L_\hi 
\right)\log\left(\frac{\Psi^0}{\epsilon}\right)\right).
\end{equation}
\end{corollary}

\begin{corollary}[Uniform minibatch sampling, choosing $\gamma$]\label{cor1usf}
In the conditions of Corollary \ref{cor1us}, we can choose $\gamma$ to obtain the best complexity: with 
\begin{equation}
\gamma = \min\left(\frac{1}{L_\ff}, \sqrt{\frac{s}{n (\max_i L_\hi) (\mu_\ff+\mu_\g+\min_i \mu_\hi)}}\right)
\end{equation}
(with the first term ignored if $L_\ff =0$),
 the iteration complexity of \algno is 
\begin{equation}
\mathcal{O}\left(\left(\frac{L_\ff}{\mu_\ff + \mu_\g+\min_{i\in [n]} \mu_{\hi} }+
 \sqrt{\frac{n \max_i L_\hi}{s  (\mu_\ff+\mu_\g+\min_i \mu_\hi)}}
+\frac{n}{s}
\right)\log\left(\frac{\Psi^0}{\epsilon}\right)\right).
\end{equation}
\end{corollary}

\begin{corollary}[Case $\mu_\hi\equiv 0$, choosing $\gamma$ and the $p_i$ with importance sampling]\label{cor1isa}
 In the conditions of Corollary \ref{cor1}, suppose in addition that  $\mu_{\dr{h_1}}=\cdots=\mu_{\dr{h_n}}=0$ and that 
 in \algno, $\Omega^t$ is of size $1$ for every $t\geq 0$. Then $\tilde{p}_i = p_i$ is the probability that $\Omega=\{i\}$ for every $i\in [n]$.
 Let $b_i\eqdef \max\left(1,\sqrt{\frac{L_{\hi}}{n(\mu_\ff+\mu_\g)}}\right)$ for every $i\in[n]$, and $b\eqdef \sum_{i=1}^n b_i$. Then by choosing 
 \begin{equation}
\gamma\eqdef\min\left(\frac{1}{L_\ff},\frac{1}{b} \max\left(\sqrt{\frac{n}{(\max_{i\in[n]} L_\hi)(\mu_\ff+\mu_\g)}},\frac{1}{\mu_\ff+\mu_\g}\right) 
\right)
\end{equation}
 (with the first term ignored if $L_\ff =0$)  and 
 $p_i\eqdef \frac{b_i}{b}$  for every $i\in[n]$, the iteration complexity of \algno is
\begin{align}
\mathcal{O}\left(\left(\frac{L_\ff}{\mu_\ff + \mu_\g}+n+\sqrt{\frac{n\bar{L}_\h}{\mu_\ff+\mu_\g}}\right)\log \left(\frac{\Psi^0}{\epsilon}\right)\right),
\end{align}
where $\bar{L}_\h=\left(\frac{1}{n}\sum_{i=1}^n \sqrt{L_\hi}\right)^2$. 
This is better than uniform sampling with $p_1=\cdots=p_n=\frac{1}{n}$ as in Corollary \ref{cor1usf} with $s=1$, since the complexity now depends on $\bar{L}_\h$ instead of $\max_i L_\hi$. 
 \end{corollary}

We note that in the conditions of Corollary \ref{cor1isa}, \algno reverts to \algn{SDM} proposed in \citet{mis19}. However, there is a mistake in their Lemma 6, so that the claim after their Theorem 6 that the $p_i$ should be proportional to the $L_\hi$ is incorrect.\bigskip

If $p_\emptyset>0$ and we want to obtain linear convergence from the $\mu_\hi>0$ but these values are unknown, we have to choose $\hat{p}=\bar{p}=0$. This case is presented next.

\begin{corollary}[Case $\hat{p}=0$]\label{cor2}
Suppose that $L_\hi<+\infty$ for every $i\in[n]$. 
In \algno, suppose that $\gamma_t\equiv \gamma$ for some $0<\gamma< \frac{2}{L_{\ff}}$ (or just $\gamma>0$ if $\ff=0$) and $\eta_i= \frac{1}{n \tilde{p}_i }$ for every $i\in[n]$.
Let \begin{equation}
\hat{\mu}_\h\eqdef \min_{i\in [n]}  \frac{2  \mu_{\hi} L_\hi}{n \tilde{p}_i(L_\hi+\mu_{\hi})}.
\end{equation}
Suppose that $\mu_{\ff}>0$ or $\mu_{\g}>0$ or $\hat{\mu}_\h>0$. Then the solution $x^\star$ of \eqref{eqpbm} exists and is unique. 
Define the Lyapunov function, for every $t\geq 0$,
\begin{equation}
\Psi^t\eqdef(1+\gamma \hat{\mu}_\h)\sqn{x^{t}-x^\star} +\frac{1-p_\emptyset}{n}\sum_{i=1}^n \frac{1}{p_i}
\left(\gamma^2  \eta_i + 
\frac{2\gamma  }{L_{\hi}+\mu_{\hi} } \right)
\sqn{u_i^{t}- \nabla \hi(x^\star)}.
\end{equation}
Then \algno converges linearly:  for every $t\geq 0$, 
$\Exp{\Psi^{t}}\leq \rho^t \Psi^0$, 
where 
\begin{align}
\rho&\eqdef\max\left(p_\emptyset+(1-p_\emptyset)\frac{\max(1-\gamma\mu_{\ff},\gamma L_{\ff}-1)^2}{(1+\gamma\mu_{\g})(1+\gamma \hat{\mu}_\h)},
1-\min_{i\in[n]}  
\frac{2n \tilde{p}_i p_i}{\gamma(L_{\hi}+\mu_{\hi})+2n \tilde{p}_i}
\right)<1.
\end{align}
Also, $x^t$ converges  to $x^\star$ and  $u_i^t$ converges to $\nabla \hi(x^\star)$ for every $i\in [n]$, 
almost surely.

As a consequence, the iteration complexity of \algno is 
\begin{equation}
\mathcal{O}\left(\left(
\frac{1}{(1-p_\emptyset)\gamma(\mu_\ff + \mu_\g+\min_{i\in [n]}  \frac{\mu_{\hi}}{n\tilde{p}_i}) }
+\max_{i\in[n]} \left(\frac{1}{p_i}+\frac{\gamma  L_\hi}{n\tilde{p}_i p_i}\right)
\right)\log\left(\frac{\Psi^0}{\epsilon}\right)\right).
\end{equation}
\end{corollary}

\begin{corollary}[Case $n=1$]\label{cor8}
Suppose that $n=1$, so that $\tilde{p}_1=1$ and $p_1=1-p_\emptyset$, and that $L_{\dr{h_1}}<+\infty$. 
Suppose that $\mu_{\ff}>0$ or $\mu_{\g}>0$ or $\mu_{\dr{h_1}}>0$. Then the solution $x^\star$ of \eqref{eqpbm} exists and is unique. 
In \algno, suppose that $\gamma_t\equiv \gamma$ for some $0<\gamma< \frac{2}{L_{\ff}}$ (or just $\gamma>0$ if $\ff=0$), $\hat{p}= \frac{1-p_\emptyset}{1-p_\emptyset+\gamma  \mu_{\dr{h_1}}}$ and 
$\eta_1= \frac{1}{1-p_\emptyset}$.
Define the Lyapunov function, for every $t\geq 0$,
\begin{equation}
\Psi^t\eqdef\frac{1-p_\emptyset+\gamma  \mu_{\dr{h_1}}}{1-p_\emptyset}\sqn{x^{t}-x^\star} +
\frac{1}{1-p_\emptyset}\left(\frac{\gamma^2}{1-p_\emptyset}  + 
\frac{2\gamma  }{L_{\dr{h_1}}+\mu_{\dr{h_1}} } \right)
\sqn{u_1^{t}- \nabla {\dr{h_1}}(x^\star)}.
\end{equation}
Then \algno converges linearly:  for every $t\geq 0$, 
$\Exp{\Psi^{t}}\leq \rho^t \Psi^0$, 
where 
\begin{align}
\rho&\eqdef\max\left(\frac{p_\emptyset\gamma \mu_{\dr{h_1}}}{1-p_\emptyset+\gamma \mu_{\dr{h_1}}}+
\frac{(1-p_\emptyset)\max(1-\gamma\mu_{\ff},\gamma L_{\ff}-1)^2}{(1+\gamma\mu_{\g})(1-p_\emptyset+\gamma \mu_{\dr{h_1}})},
1-
\frac{2(1-p_\emptyset)^2}{\gamma(L_{\dr{h_1}}+\mu_{\dr{h_1}})+2(1-p_\emptyset)}
\right)<1.\label{eqrho6}
\end{align}
As a consequence, the iteration complexity of \algno is 
\begin{equation}
\mathcal{O}\left(\left(
\frac{1}{\frac{(1-p_\emptyset )(1+\gamma \mu_{\dr{h_1}})}{1-p_\emptyset+\gamma \mu_{\dr{h_1}}}\gamma(\mu_\ff + \mu_\g)+\gamma \mu_{\dr{h_1}}  }
+\frac{1-p_\emptyset+\gamma L_{\dr{h_1}}}{(1-p_\emptyset)^2}
\right)\log\left(\frac{\Psi^0}{\epsilon}\right)\right).
\end{equation}
\end{corollary}

In the conditions of Corollary \ref{cor8}, \algno reverts to an instance of \algn{RandProx-Skip} proposed in \citet{con23rp}. If $\mu_{\dr{h_1}}=0$, we recover the exact same rate as in \citet[Theorem 1]{con23rp}. Thus, our analysis is more general as it exploits  strong convexity of ${\dr{h_1}}$.
We also show in Section \ref{seccdy} that our results are comparable to existing ones for the deterministic \algn{Davis--Yin algorithm}.

\section{Proof of Corollary \ref{cor1isa} }

In the conditions of Corollary \ref{cor1} with the additional assumptions in Corollary \ref{cor1isa}, 
 the iteration complexity of \algno is 
\begin{equation*}
\mathcal{O}\left(\left(
\frac{1}{\gamma  (\mu_\ff+\mu_\g) }
+\max_{i\in[n]} \left(\frac{1}{p_i}+\frac{\gamma  L_\hi}{n p_i^2}\right)
\right)\log\left(\frac{\Psi^0}{\epsilon}\right)\right).
\end{equation*}
By choosing 
\begin{equation*}
\gamma=\min\left(\frac{1}{L_\ff},\min_{i\in[n]} \frac{p_i\sqrt{n}}{\sqrt{L_\hi(\mu_\ff+\mu_\g)}}\right)
\end{equation*}
(with the first term ignored if $L_\ff =0$), 
 the complexity becomes
\begin{equation*}
\mathcal{O}\left(\left(\frac{L_\ff}{\mu_\ff + \mu_\g}+\max_{i\in[n]} \frac{1}{p_i}\max\left(1,\sqrt{\frac{L_{\hi}}{n(\mu_\ff+\mu_\g)}}\right)\right)\log \left(\frac{\Psi^0}{\epsilon}\right)\right).
\end{equation*}
Let $b_i\eqdef \max\left(1,\sqrt{\frac{L_{\hi}}{n(\mu_\ff+\mu_\g)}}\right)$ for every $i\in[n]$, and $b\eqdef \sum_{i=1}^n b_i$. Then by choosing $p_i\eqdef \frac{b_i}{b}$, the complexity becomes 
\begin{align*}
\mathcal{O}\left(\left(\frac{L_\ff}{\mu_\ff + \mu_\g}+b\right)\log\left(\frac{\Psi^0}{\epsilon}\right)\right)&=\mathcal{O}\left(\left(\frac{L_\ff}{\mu_\ff + \mu_\g}+\sum_{i=1}^n b_i\right)\log \left(\frac{\Psi^0}{\epsilon}\right)\right)\\
&= \mathcal{O}\left(\left(\frac{L_\ff}{\mu_\ff + \mu_\g}+\sum_{i=1}^n \left(1+\sqrt{\frac{L_{\hi}}{n(\mu_\ff+\mu_\g)}}\right)\right)\log \left(\frac{\Psi^0}{\epsilon}\right)\right)\\
&= \mathcal{O}\left(\left(\frac{L_\ff}{\mu_\ff + \mu_\g}+n+\sum_{i=1}^n \sqrt{\frac{L_{\hi}}{n(\mu_\ff+\mu_\g)}}\right)\log \left(\frac{\Psi^0}{\epsilon}\right)\right)\\
&=\mathcal{O}\left(\left(\frac{L_\ff}{\mu_\ff + \mu_\g}+n+\sqrt{\frac{n\bar{L}_\h}{\mu_\ff+\mu_\g}}\right)\log \left(\frac{\Psi^0}{\epsilon}\right)\right),
\end{align*}
where $\bar{L}_\h=\left(\frac{1}{n}\sum_{i=1}^n \sqrt{L_\hi}\right)^2$. 
Moreover, we have 
\begin{align*}
\gamma&=\min\left(\frac{1}{L_\ff},\frac{1}{b}\min_{i\in[n]} \max\left(\sqrt{\frac{n}{L_\hi(\mu_\ff+\mu_\g)}},\frac{1}{\mu_\ff+\mu_\g}\right) 
\right)\\
&=\min\left(\frac{1}{L_\ff},\frac{1}{b} \max\left(\sqrt{\frac{n}{(\max_{i\in[n]} L_\hi)(\mu_\ff+\mu_\g)}},\frac{1}{\mu_\ff+\mu_\g}\right) 
\right).
\end{align*}

\section{Proof of Theorem \ref{theog0} }

Let $t\geq 0$. Starting from \eqref{eqp1}, since $q^t=q^\star=0$, we have
\begin{align*}
&(1+\gamma \hat{\mu}_\h)\Expc{\sqn{x^{t+1}-x^\star}} +\frac{1- p_\emptyset + \bar{p}}{1- p_\emptyset}
\left(\gamma^2  \eta_1 + 
\frac{2\gamma  }{L_{\dr{h_1}}+\mu_{\dr{h_1}}} \right)
\Expc{\sqn{u_1^{t+1}-u_1^\star}}\notag\\
&\leq    \left(p_\emptyset(1-\hat{p})+
(1- p_\emptyset + \bar{p})\frac{\max(1-\gamma\mu_{\ff},\gamma L_{\ff}-1)^2}{(1+\gamma\mu_{\g})(1+\gamma \hat{\mu}_\h)}\right)(1+\gamma \hat{\mu}_\h)\sqn{x^t-x^\star}\\
&\quad+\frac{1- p_\emptyset + \bar{p}}{1-p_\emptyset} \left(1-\frac{2(1-p_\emptyset)}{\gamma\eta_1(L_{\dr{h_1}}+\mu_{\dr{h_1}})+2}\right)
\left(\gamma^2  \eta_1 + 
\frac{2\gamma  }{L_{\dr{h_1}}+\mu_{\dr{h_1}}} \right)\sqn{u_1^{t} -u_1^\star }\notag\\
&\quad- \frac{\gamma^2(1- p_\emptyset + \bar{p})}{1+\gamma\mu_{\g}}\sqn{ u_1^t-u_1^\star}\\
&=   \left(p_\emptyset(1-\hat{p})+
(1- p_\emptyset + \bar{p})\frac{\max(1-\gamma\mu_{\ff},\gamma L_{\ff}-1)^2}{(1+\gamma\mu_{\g})(1+\gamma \hat{\mu}_\h)}\right)(1+\gamma \hat{\mu}_\h)\sqn{x^t-x^\star}\label{eqp2}\\
&\quad+\left(1-\frac{(1-p_\emptyset)^2\big(2(1+\gamma\mu_{\g})+\gamma(L_{\dr{h_1}}+\mu_{\dr{h_1}})\big)
}{\big((1- p_\emptyset + \bar{p})\gamma(L_{\dr{h_1}}+\mu_{\dr{h_1}})+2(1-p_\emptyset)\big)(1+\gamma\mu_{\g})}\right)\\
&\quad\times \frac{1- p_\emptyset + \bar{p}}{1-p_\emptyset} 
\left(\gamma^2  \eta_1 + 
\frac{2\gamma  }{L_{\dr{h_1}}+\mu_{\dr{h_1}}} \right)\sqn{u_1^{t} -u_1^\star}.
\end{align*}
Hence,
\begin{equation}
\Expc{\Psi^{t+1}}\leq \rho \Psi^t,\label{eqrec2bzq}
\end{equation}
with
\begin{align*}
\rho\eqdef\max&\left(p_\emptyset(1-\hat{p})+
(1- p_\emptyset + \bar{p})\frac{\max(1-\gamma\mu_{\ff},\gamma L_{\ff}-1)^2}{(1+\gamma\mu_{\g})(1+\gamma \hat{\mu}_\h)},\right.\\
&\quad\left.
1-\frac{(1-p_\emptyset)^2\big(2(1+\gamma\mu_{\g})+\gamma(L_{\dr{h_1}}+\mu_{\dr{h_1}})\big)
}{\big((1- p_\emptyset + \bar{p})\gamma(L_{\dr{h_1}}+\mu_{\dr{h_1}})+2(1-p_\emptyset)\big)(1+\gamma\mu_{\g})}
\right),
\end{align*}
where the second term is equal to 
\begin{align*}
1-\frac{(1-p_\emptyset)^2
}{(1- p_\emptyset + \bar{p})(1+\gamma\mu_{\g})}
\end{align*}
if $L_{\dr{h_1}}=+\infty$.

Using the tower rule, we can unroll the recursion in \eqref{eqrec2bzq} to obtain the unconditional expectation of $\Psi^{t+1}$.  Since $\rho<1$, $\Exp{\Psi^t}\rightarrow 0$. 
Moreover, using classical results on supermartingale convergence \citep[Proposition A.4.5]{ber15}, it follows from \eqref{eqrec2bzq} that $\Psi^t\rightarrow 0$ almost surely. Almost sure convergence of $x^t$ and the $u_i^t$ follows. 

\section{Proof of Theorem \ref{theosi1}}\label{apptheosi1}

We consider \algn{SMPM-EXT}, an extended version of \algno with additional variables $z_i^t$ defined, so that $u_i^t = \nabla \hi(z_{i}^{t+1})$ for every $i\in [n]$ and $t\geq 0$. The variables $z_i^t$ are not actually computed or stored, they are just defined for the analysis.
 
 \begin{figure*}[!t]	
\begin{algorithm}[H]
		\caption{\algn{SMPM-EXT}}
		\begin{algorithmic}[1]
			\STATE  \textbf{parameters:} initial points $x^0,z_1^0,\ldots,z_n^0\in\mathcal{X}$; 
			probability $\hat{p}\in[0,1]$
			\STATE 
			sampling distribution $\mathcal{S}$; 
			stepsizes $(\gamma_t)_{t\geq 0}$, $(\eta_i)_{i=1}^n$
			\STATE  let $u_i^0 \in \partial  \hi(z_{i}^{0})$, $\forall i\in[n]$
			\STATE $\vv^0\coloneqq\frac{1}{n} \sum_{i=1}^n u_i^0$
			\FOR{$t=0, 1, \ldots$}
			\STATE $\hat{x}^{t} \coloneqq  \mathrm{prox}_{\gamma_t \g}\big(x^t -\gamma_t \nabla \ff(x^t) - \gamma_t \vv^t\big)$
			\STATE sample $\Omega^t\sim \mathcal{S}$ 
			\IF{$\Omega^t\neq \emptyset$}
			\FOR{$i\in\Omega^t$}
			\STATE $y_{i}^{t+1}\coloneqq \mathrm{prox}_{\gamma_t  \eta_i\hi} (\hat{x}^{t} + \gamma_t  \eta_i u_{i}^t )$
			\STATE $z_i^{t+1}\eqdef y_i^{t+1}$
			\STATE $u_{i}^{t+1}\coloneqq u_{i}^t+\frac{1}{\gamma_t \eta_i}  \big(\hat{x}^{t}- y_{i}^{t+1}\big)\in \partial  \hi(y_{i}^{t+1})$ 
			\ENDFOR
			\STATE $u_{i}^{t+1}\coloneqq u_{i}^t$, $z_i^{t+1}\eqdef z_i^t$, $\forall i\in[n]\backslash \Omega^t$ 
			\STATE  $x^{t+1} \coloneqq \frac{1}{|\Omega^t|}\sum_{i \in \Omega^t}  y_{i}^{t+1}$ 
			\STATE  $\vv^{t+1}\coloneqq 
			\frac{1}{n}\sum_{i=1}^nu_i^{t+1}$
			\ELSE
			\STATE $x^{t+1}\coloneqq \left\{\begin{array}{l} \hat{x}^t\ \mbox{with probability}\ \hat{p}\\
			x^t\ \mbox{with probability}\ 1-\hat{p}
			\end{array}\right.$
			\STATE $u_{i}^{t+1}\coloneqq u_{i}^t$ and $z_i^{t+1}\eqdef z_i^t$, $\forall i\in[n]\backslash \Omega^t$ 
			\ENDIF
			\ENDFOR
		\end{algorithmic}
		\end{algorithm}\end{figure*}

Let $t\geq 0$. Following the derivations in the beginning of Section \ref{secap1}, we have from \eqref{eqy7}, \eqref{eqco4} and \eqref{eqw1}, 
\begin{align*}
&(1+\gamma\hat{\mu}_\h)\frac{1}{n}\sum_{i=1}^n \sqn{y_{i}^{t+1}-x^\star} +\left(\gamma^2   +
\frac{2 \gamma  }{L_{\h}+\mu_{\h} } \right)\frac{1}{n p_s}\sum_{i=1}^n 
\Expc{\sqn{u_i^{t+1}-u_i^\star}}\notag\\
&\leq \max(1-\gamma\mu_{\ff},\gamma L_{\ff}-1)^2 \sqn{x^t-x^\star} + \frac{\gamma^2}{n}\sum_{i=1}^n \sqn{u_i^{t} -u_i^\star }\\
&\quad+\left(\gamma^2   +
\frac{2 \gamma  }{L_{\h}+\mu_{\h} } \right)\frac{1-p_s}{n p_s}\sum_{i=1}^n \sqn{u_i^{t} -u_i^\star }- \gamma^2\sqn{\frac{1}{n} \sum_{j=1}^n (u_j^t-u_j^\star)}.
\end{align*}
Let $\alpha \in [0,1]$.		
We then have 
\begin{align}
&(1+\gamma \hat{\mu}_\h)\frac{1}{n}\sum_{i=1}^n \sqn{y_{i}^{t+1}-x^\star} +\left(\gamma^2   + 
\frac{2\gamma  }{L_{\h}+\mu_{\h} } \right)\frac{1}{n p_s}\sum_{i=1}^n
\Expc{\sqn{u_i^{t+1}-u_i^\star}}\notag\\
&\leq 
\max(1-\gamma\mu_{\ff},\gamma L_{\ff}-1)^2\sqn{x^t-x^\star}\notag\\
&\quad+\left(\left(
\gamma^2+\frac{2\gamma}{L_{\h}+\mu_{\h}}
\right)\frac{1-p_s}{p_s}+(1-\alpha)\gamma^2\right)\frac{1}{n}\sum_{i=1}^n  \sqn{u_i^{t} -u_i^\star }\notag\\
&+\frac{\alpha\gamma^2}{n}\sum_{i=1}^n \sqn{u_i^{t} -u_i^\star -\frac{1}{n} \sum_{j=1}^n (u_j^t-u_j^\star)}\notag\\
&=
\max(1-\gamma\mu_{\ff},\gamma L_{\ff}-1)^2\sqn{x^t-x^\star}\notag\\
&\quad+\left(\left(
\gamma^2+\frac{2\gamma}{L_{\h}+\mu_{\h}}
\right)\frac{1-p_s}{p_s}+(1-\alpha)\gamma^2\right)\frac{1}{n}\sum_{i=1}^n  \sqn{u_i^{t} -u_i^\star }\notag\\
&+\frac{\alpha\gamma^2}{n}\sum_{i=1}^n \sqn{\left(\nabla \hi(z_i^{t})-\frac{1}{n} \sum_{j=1}^n \nabla \dr{h_j}(z_j^{t})\right)-\left(
 \nabla \hi(x^\star) -\frac{1}{n} \sum_{j=1}^n \nabla \dr{h_j}(x^{\star})\right)}\notag\\
 &\leq 
\max(1-\gamma\mu_{\ff},\gamma L_{\ff}-1)^2\sqn{x^t-x^\star}+ \frac{\delta^2\alpha\gamma^2}{n}\sum_{i=1}^n \sqn{z_i^t-x^\star}\notag\\
&\quad+\left(\left(
\gamma^2+\frac{2\gamma}{L_{\h}+\mu_{\h}}
\right)\frac{1-p_s}{p_s}+(1-\alpha)\gamma^2\right)\frac{1}{n}\sum_{i=1}^n\sqn{u_i^{t} -u_i^\star }.\label{eqsad}
\end{align}
For every $i\in[n]$, we have
\begin{equation*}
\Expc{\sqn{z_i^{t+1}-x^\star}}=p_s\sqn{y_i^{t+1}-x^\star}+(1-p_s)\sqn{z_i^{t}-x^\star},
\end{equation*}
so that
\begin{equation}
\frac{1}{n p_s}\sum_{i=1}^n 
\Expc{\sqn{z_i^{t+1}-x^\star}}-\frac{1}{n}\sum_{i=1}^n \sqn{y_i^{t+1}-x^\star}= \frac{1-p_s}{n p_s}\sum_{i=1}^n \sqn{z_i^{t}-x^\star}.\label{eqrro}
\end{equation}
By adding $\frac{\gamma}{2}(\hat{\mu}_\ff+\hat{\mu}_\h)$ times \eqref{eqrro} to \eqref{eqsad}, we obtain
\begin{align*}
&\left(1-\frac{\gamma\hat{\mu}_\ff}{2}+\frac{\gamma\hat{\mu}_\h}{2} \right)\frac{1}{n}\sum_{i=1}^n \sqn{y_{i}^{t+1}-x^\star}
+\frac{\gamma(\hat{\mu}_\ff+\hat{\mu}_\h)}{2n p_s}\sum_{i=1}^n 
\Expc{\sqn{z_i^{t+1}-x^\star}}\\
& \quad+\left(\gamma^2   + 
\frac{2\gamma  }{L_{\h}+\mu_{\h} } \right)\frac{1}{n p_s}\sum_{i=1}^n
\Expc{\sqn{u_i^{t+1}-u_i^\star}}\notag\\
&\leq 
(1-\gamma\hat{\mu}_\ff)\sqn{x^t-x^\star}+ \left(\frac{\delta^2\alpha\gamma^2}{n}+\frac{\gamma(\hat{\mu}_\ff+\hat{\mu}_\h)(1-p_s)}{2n p_s}\right)
\sum_{i=1}^n \sqn{z_i^t-x^\star}\notag\\
&\quad+\left(\left(
\gamma^2+\frac{2\gamma}{L_{\h}+\mu_{\h}}
\right)\frac{1-p_s}{p_s}+(1-\alpha)\gamma^2\right)\frac{1}{n}\sum_{i=1}^n\sqn{u_i^{t} -u_i^\star }\\
&=
(1-\gamma\hat{\mu}_\ff)\sqn{x^t-x^\star}+ \left(1-p_s+\frac{2p_s\delta^2\alpha\gamma}{\hat{\mu}_\ff+\hat{\mu}_\h}\right)\frac{\gamma(\hat{\mu}_\ff+\hat{\mu}_\h)}{2n p_s}
\sum_{i=1}^n \sqn{z_i^t-x^\star}\notag\\
&\quad+\left(1-p_s+(1-\alpha)p_s
\frac{\gamma(L_\h+\mu_\h)}{\gamma(L_\h+\mu_\h)+2}
\right)\left(
\gamma^2+\frac{2\gamma}{L_{\h}+\mu_{\h}}
\right)\frac{1}{np_s}\sum_{i=1}^n\sqn{u_i^{t} -u_i^\star }.
\end{align*}
We choose $\alpha$ to balance the two contraction factors in the $z$ and $u$ terms; that is, to have
\begin{align*}
&\alpha\frac{2\delta^2\gamma}{\hat{\mu}_\ff+\hat{\mu}_\h}=(1-\alpha)
\frac{\gamma(L_\h+\mu_\h)}{\gamma(L_\h+\mu_\h)+2}\\
\Leftrightarrow\ &
\alpha=\left(1+\frac{2\delta^2\gamma}{\hat{\mu}_\ff+\hat{\mu}_\h}\frac{\gamma(L_\h+\mu_\h)+2}{\gamma(L_\h+\mu_\h)}\right)^{-1}=\frac{(\hat{\mu}_\ff+\hat{\mu}_\h)(L_\h+\mu_\h)}{(\hat{\mu}_\ff+\hat{\mu}_\h)(L_\h+\mu_\h)+2\delta^2\big(\gamma(L_\h+\mu_\h)+2\big)}.
\end{align*}
With this value of $\alpha$, and using
\begin{align*}
\Expc{\sqn{x^{t+1}-x^\star}}&\leq \frac{1}{n}\sum_{i=1}^n \sqn{y_{i}^{t+1}-x^\star}.
\end{align*}
shown in \eqref{eqco2}, 
we obtain
\begin{align*}
&\left(1-\frac{\gamma\hat{\mu}_\ff}{2}+\frac{\gamma\hat{\mu}_\h}{2} \right)\Expc{\sqn{x^{t+1}-x^\star}}
+\frac{\gamma(\hat{\mu}_\ff+\hat{\mu}_\h)}{2n p_s}\sum_{i=1}^n 
\Expc{\sqn{z_i^{t+1}-x^\star}}\\
& \quad+\left(\gamma^2   + 
\frac{2\gamma  }{L_{\h}+\mu_{\h} } \right)\frac{1}{n p_s}\sum_{i=1}^n
\Expc{\sqn{u_i^{t+1}-u_i^\star}}\notag\\
&\leq 
\frac{2-2\gamma\hat{\mu}_\ff}{2-\gamma\hat{\mu}_\ff+\gamma\hat{\mu}_\h}\left(1-\frac{\gamma\hat{\mu}_\ff}{2}+\frac{\gamma\hat{\mu}_\h}{2} \right)\sqn{x^t-x^\star}\\
&\quad+\left(1-p_s\frac{(\hat{\mu}_\ff+\hat{\mu}_\h)(L_\h+\mu_\h)+4\delta^2}{(\hat{\mu}_\ff+\hat{\mu}_\h)(L_\h+\mu_\h)+4\delta^2+2\delta^2\gamma(L_\h+\mu_\h)}
\right)\frac{\gamma(\hat{\mu}_\ff+\hat{\mu}_\h)}{2n p_s}
\sum_{i=1}^n \sqn{z_i^t-x^\star}\notag\\
&\quad+\left(1-p_s\frac{(\hat{\mu}_\ff+\hat{\mu}_\h)(L_\h+\mu_\h)+4\delta^2}{(\hat{\mu}_\ff+\hat{\mu}_\h)(L_\h+\mu_\h)+4\delta^2+2\delta^2\gamma(L_\h+\mu_\h)}
\right)\left(
\gamma^2+\frac{2\gamma}{L_{\h}+\mu_{\h}}
\right)\frac{1}{np_s}\sum_{i=1}^n\sqn{u_i^{t} -u_i^\star }.
\end{align*}
Hence
\begin{equation*}
\Expc{\Psi^{t+1}}\leq \max\left(\frac{2-2\gamma\hat{\mu}_\ff}{2-\gamma\hat{\mu}_\ff+\gamma\hat{\mu}_\h},
1-p_s\frac{(\hat{\mu}_\ff+\hat{\mu}_\h)(L_\h+\mu_\h)+4\delta^2}{(\hat{\mu}_\ff+\hat{\mu}_\h)(L_\h+\mu_\h)+4\delta^2+2\delta^2\gamma(L_\h+\mu_\h)}\right)\Psi^t.
\end{equation*}
Unrolling the recursion, we obtain the unconditional expectation of $\Psi^t$.

\section{Proof of Theorem~\ref{theoac}}

Let $t\geq 0$. Following the derivations in the beginning of Section \ref{secap1}, we have in \eqref{eqp1}
\begin{align*}
\Expc{\Psi^{t+1}}
&=(1+\gamma_t \hat{\mu}_\h)\Expc{\sqn{x^{t+1}-x^\star}} +\frac{1}{n}\sum_{i=1}^n \frac{\gamma_t^2  \eta_i }{p_i}
\Expc{\sqn{u_i^{t+1}-u_i^\star}}\notag\\
&\leq   
\frac{(1-\gamma_t\mu_{\ff})^2}{(1+\gamma_t\mu_{\g})(1+\gamma_t \hat{\mu}_\h)}(1+\gamma_t \hat{\mu}_\h)\sqn{x^t-x^\star}+\frac{1}{n}\sum_{i=1}^n
\frac{\gamma_t^2  \eta_i }{p_i}
\sqn{u_i^{t} -u_i^\star }\\
&=
C_t(1+\gamma_{t-1} \hat{\mu}_\h)\sqn{x^t-x^\star}+\frac{\gamma_t^2}{\gamma_{t-1}^2}\frac{1}{n}\sum_{i=1}^n
\frac{\gamma_{t-1}^2  \eta_i }{p_i}
\sqn{u_i^{t} -u_i^\star }
\end{align*}
with
\begin{equation}
C_t\eqdef\frac{(1-\gamma_t\mu_{\ff})^2}{(1+\gamma_t\mu_{\g})(1+\gamma_{t-1} \hat{\mu}_\h)}.
\end{equation}
Moreover, $C_t\leq 1-\gamma_t \mu$, because  $\gamma_{t-1}\geq \gamma_t$ and
$\frac{1}{1+z}\leq 1-\frac{z}{2}$ for any $z\in[0,1]$, and we have
\begin{align}
\left(1-\gamma_t \mu\right)\frac{\gamma_{t-1}^2}{\gamma_t^2}=\frac{\big(1-\frac{2}{a+t}\big)(a+t)^2}{(a+t-1)^2}=\frac{(a+ t-2)(a+ t)}{(a+t-1)^2}<1,
\end{align}
so that $C_t\leq \frac{\gamma_{t}^2}{\gamma_{t-1}^2}$ and
\begin{align*}
\Expc{\Psi^{t+1}}&\leq \frac{\gamma_t^2}{\gamma_{t-1}^2}\Psi^t.
\end{align*}
Unrolling the recursion, we obtain, for every $t\geq 0$,
\begin{align*}
\Exp{\Psi^{t}}&\leq \frac{\gamma_{t-1}^2}{\gamma_{-1}^2}\Psi^0.
\end{align*}

\section{Proof of Theorem~\ref{theoac3}}

Let $t\geq 0$. Following the derivations in the beginning of Section \ref{secap1}, we have in \eqref{eqp1}
\begin{align*}
\Expc{\Psi^{t+1}}
&=\Expc{\sqn{x^{t+1}-x^\star}} +\frac{1}{n}\sum_{i=1}^n \frac{\gamma_t^2  \eta_i }{p_i}
\Expc{\sqn{u_i^{t+1}-u_i^\star}}\notag\\
&\leq   
\frac{(1-\gamma_t\mu_{\ff})^2}{1+\gamma_t\mu_{\g}}\sqn{x^t-x^\star}+\frac{1}{n}\sum_{i=1}^n
\frac{\gamma_t^2  \eta_i }{p_i}
\sqn{u_i^{t} -u_i^\star }\\
&=
C_t\sqn{x^t-x^\star}+\frac{\gamma_t^2}{\gamma_{t-1}^2}\frac{1}{n}\sum_{i=1}^n
\frac{\gamma_{t-1}^2  \eta_i }{p_i}
\sqn{u_i^{t} -u_i^\star }
\end{align*}
with
\begin{equation}
C_t\eqdef\frac{(1-\gamma_t\mu_{\ff})^2}{1+\gamma_t\mu_{\g}}.
\end{equation}
Moreover, $C_t\leq 1-\gamma_t \mu$, because 
$\frac{1}{1+z}\leq 1-\frac{z}{2}$ for any $z\in[0,1]$, and we have
\begin{align}
\left(1-\gamma_t \mu\right)\frac{\gamma_{t-1}^2}{\gamma_t^2}=\frac{\big(1-\frac{2}{a+t}\big)(a+t)^2}{(a+t-1)^2}=\frac{(a+ t-2)(a+ t)}{(a+t-1)^2}<1,
\end{align}
so that $C_t\leq \frac{\gamma_{t}^2}{\gamma_{t-1}^2}$ and
\begin{align*}
\Expc{\Psi^{t+1}}&\leq \frac{\gamma_t^2}{\gamma_{t-1}^2}\Psi^t.
\end{align*}
Unrolling the recursion, we obtain, for every $t\geq 0$,
\begin{align*}
\Exp{\Psi^{t}}&\leq \frac{\gamma_{t-1}^2}{\gamma_{-1}^2}\Psi^0.
\end{align*}

\section{$\mathcal{O}(1/t)$ Convergence in the General Convex Case} 

\begin{theorem}[$\mathcal{O}(1/t)$ convergence]\label{theogen}
Suppose that in \algno, $\hat{p}\in [0,1]$, $\gamma_t\equiv \gamma$ for some $0<\gamma< \frac{2}{L_{\ff}}$ (or just $\gamma>0$ if $\ff=0$),
and  for every $i\in[n]$, $\eta_i= \frac{1- p_\emptyset + \bar{p}}{n \tilde{p}_i (1-p_\emptyset)}$.
Let $x^\star$ and $u_i^\star \in \partial\hi(x^\star)$, $i=1,\ldots,n$, satisfy \eqref{eqincl}. 
Let $T\geq 0$ and choose $t$ uniformly at random in $\{0,\ldots,T\}$. Then
\begin{align*}
&\gamma(2-\gamma L_\ff)(1- p_\emptyset + p_\emptyset\hat{p})\Exp{\langle x^t-x^\star,\nabla\ff(x^t)-\nabla\ff(x^\star)\rangle}\\
&+(1\!-\! p_\emptyset \!+\! p_\emptyset\hat{p})\gamma^2\Exp{\sqn{\tilde{\nabla}\g(\hat{x}^{t})\!-\!\tilde{\nabla}\g(\hat{x}^\star)+\frac{1}{n} \sum_{j=1}^n (u_j^t\!-\!u_j^\star)}}\\
&+2\gamma(1- p_\emptyset + p_\emptyset\hat{p})\Exp{ \langle \tilde{\nabla}\g(\hat{x}^{t})-\tilde{\nabla}\g(\hat{x}^\star), \hat{x}^{t} -x^\star\rangle} \\
&+2\gamma \frac{1- p_\emptyset +p_\emptyset\hat{p}}{n}\sum_{i=1}^n  \Exp{ \langle \tilde{\nabla} \hi(y_i^{t+1}) - u_i^\star,y_{i}^{t+1}-x^\star\rangle}\\
&\leq \frac{1}{T+1}\left(\sqn{x^0\!-\!x^\star} +
\frac{1\!-\! p_\emptyset\! +\! p_\emptyset\hat{p}}{n}\gamma^2 \sum_{i=1}^n 
\frac{ \eta_i}{p_i}
\sqn{u_i^{0}\! -\!u_i^\star }\right).
\end{align*}
\end{theorem}

Proof: Let $t\geq 0$. Starting with the same derivations as in Section \ref{secap1}, from  \eqref{eqco1}, we have, for every $i\in[n]$,
\begin{align*}
&\sqn{y_{i}^{t+1}-x^\star}+\gamma^2  \eta_i^2\left(\frac{1}{p_i}\Expc{\sqn{u_i^{t+1}-u_i^\star}}-\frac{1-p_i}{p_i}\sqn{u_i^{t} -u_i^\star }\right)\\
&=\sqn{y_{i}^{t+1}-x^\star}+\gamma^2  \eta_i^2 \sqn{\tilde{\nabla} \hi(y_i^{t+1}) - u_i^\star} \notag\\
&= \sqn{\hat{x}^{t} -x^\star+ \gamma  \eta_i (u_{i}^t-u_i^\star)}-2\gamma  \eta_i \langle
y_{i}^{t+1}-x^\star,\tilde{\nabla} \hi(y_i^{t+1}) - u_i^\star\rangle.\notag
\end{align*}
Combining this equality  
with \eqref{eqco2}, we obtain
\begin{align*}
&\frac{1}{1- p_\emptyset}\left(\Expc{\sqn{x^{t+1}-x^\star}}- p_\emptyset\hat{p}\sqn{\hat{x}^t-x^\star}-p_\emptyset(1-\hat{p})\sqn{x^t-x^\star}\right)\\
&\quad +\sum_{i=1}^n  \tilde{p}_i\gamma^2  \eta_i^2 
\left(\frac{1}{p_i}\Expc{\sqn{u_i^{t+1}-u_i^\star}}-\frac{1-p_i}{p_i}\sqn{u_i^{t} -u_i^\star }\right)\\
&\leq \sum_{i=1}^n \tilde{p}_i\sqn{y_{i}^{t+1}-x^\star} +\sum_{i=1}^n  \tilde{p}_i\gamma^2  \eta_i^2
\left(\frac{1}{p_i}\Expc{\sqn{u_i^{t+1}-u_i^\star}}-\frac{1-p_i}{p_i}\sqn{u_i^{t} -u_i^\star }\right)\\
&\leq  \sum_{i=1}^n\tilde{p}_i \sqn{\hat{x}^{t} -x^\star+ \gamma  \eta_i (u_{i}^t-u_i^\star)}-2\gamma   \sum_{i=1}^n\tilde{p}_i \eta_i \langle
y_{i}^{t+1}-x^\star,\tilde{\nabla} \hi(y_i^{t+1}) - u_i^\star\rangle,
\end{align*}
so that
\begin{align}
&\frac{1}{1- p_\emptyset}\Expc{\sqn{x^{t+1}-x^\star}} +\sum_{i=1}^n \frac{\tilde{p}_i}{p_i}
\gamma^2  \eta_i^2 
\Expc{\sqn{u_i^{t+1}-u_i^\star}}\notag\\
&\leq \frac{p_\emptyset \hat{p}}{1- p_\emptyset} \sqn{\hat{x}^{t}-x^\star}+\frac{p_\emptyset(1-\hat{p})}{1-p_\emptyset}\sqn{x^t-x^\star}+\sum_{i=1}^n 
\frac{ \tilde{p}_i(1-p_i)}{p_i}\gamma^2  \eta_i^2 
\sqn{u_i^{t} -u_i^\star }\notag\\
&\quad+ \sum_{i=1}^n \tilde{p}_i\sqn{\hat{x}^{t} -x^\star+ \gamma  \eta_i (u_{i}^t-u_i^\star)}-2\gamma  \sum_{i=1}^n\tilde{p}_i  \eta_i \langle
y_{i}^{t+1}-x^\star,\tilde{\nabla} \hi(y_i^{t+1}) - u_i^\star\rangle.\label{eqpp0p}
\end{align}
We have supposed that for every $i\in[n]$,
\begin{equation*}
\eta_i= \frac{1- p_\emptyset + p_\emptyset\hat{p}}{n \tilde{p}_i (1-p_\emptyset)}.
\end{equation*}
So, multiplying \eqref{eqpp0p} by $1- p_\emptyset$, we obtain
\begin{align*}
&\Expc{\sqn{x^{t+1}-x^\star}} +\frac{1- p_\emptyset +p_\emptyset\hat{p}}{n}\sum_{i=1}^n \frac{1}{p_i}
\gamma^2  \eta_i 
\Expc{\sqn{u_i^{t+1}-u_i^\star}}\notag\\
&\leq p_\emptyset(1-\hat{p})\sqn{x^t-x^\star}+
\frac{1- p_\emptyset + p_\emptyset\hat{p}}{n}\sum_{i=1}^n 
\frac{ 1-p_i}{p_i}\gamma^2  \eta_i 
\sqn{u_i^{t} -u_i^\star }\\
&\quad + p_\emptyset\hat{p} \sqn{\hat{x}^{t}-x^\star}+ (1- p_\emptyset)\sum_{i=1}^n \tilde{p}_i\sqn{\hat{x}^{t} -x^\star+ \gamma  \eta_i (u_{i}^t-u_i^\star)}\\
&\quad-2\gamma \frac{1- p_\emptyset +p_\emptyset\hat{p}}{n}\sum_{i=1}^n   \langle y_{i}^{t+1}-x^\star,\tilde{\nabla} \hi(y_i^{t+1}) - u_i^\star\rangle\\
&= p_\emptyset(1-\hat{p})\sqn{x^t-x^\star}+
\frac{1- p_\emptyset +p_\emptyset\hat{p}}{n}\sum_{i=1}^n 
\frac{ 1-p_i}{p_i}\gamma^2  \eta_i  
\sqn{u_i^{t} -u_i^\star }\\
&\quad + (1- p_\emptyset + p_\emptyset\hat{p})\left(\frac{p_\emptyset\hat{p}}{1- p_\emptyset + p_\emptyset\hat{p}} \sqn{\hat{x}^{t}-x^\star}+ \frac{1- p_\emptyset}{1- p_\emptyset + p_\emptyset\hat{p}}\sum_{i=1}^n \tilde{p}_i\sqn{\hat{x}^{t} -x^\star+ \gamma  \eta_i (u_{i}^t-u_i^\star)}\right)\\
&\quad-2\gamma \frac{1- p_\emptyset +p_\emptyset\hat{p}}{n}\sum_{i=1}^n   \langle y_{i}^{t+1}-x^\star,\tilde{\nabla} \hi(y_i^{t+1}) - u_i^\star\rangle.
\end{align*}
We have shown in \eqref{eqco4} that
\begin{align*}
&\frac{p_\emptyset\hat{p}}{1- p_\emptyset + p_\emptyset\hat{p}} \sqn{\hat{x}^{t}-x^\star}+ \frac{1- p_\emptyset}{1- p_\emptyset + p_\emptyset\hat{p}}\sum_{i=1}^n \tilde{p}_i\sqn{\hat{x}^{t} -x^\star+ \gamma  \eta_i (u_{i}^t-u_i^\star)}\\
&=\sqn{w^t-w^\star}- \gamma^2\sqn{\tilde{\nabla}\g(\hat{x}^{t})-\tilde{\nabla}\g(\hat{x}^\star)+\frac{1}{n} \sum_{j=1}^n (u_j^t-u_j^\star)}
\\
&\quad- 2\gamma \langle \tilde{\nabla}\g(\hat{x}^{t})-\tilde{\nabla}\g(\hat{x}^\star), \hat{x}^{t} -x^\star\rangle  +\frac{ \gamma^2}{n} \sum_{i=1}^n   \eta_i \sqn{u_{i}^t-u_i^\star}.
 \end{align*}
Therefore,
\begin{align*}
&\Expc{\sqn{x^{t+1}-x^\star}} +\frac{1- p_\emptyset + p_\emptyset\hat{p}}{n}\sum_{i=1}^n \frac{1}{p_i}
\gamma^2  \eta_i 
\Expc{\sqn{u_i^{t+1}-u_i^\star}}\notag\\
&\leq  p_\emptyset(1-\hat{p})\sqn{x^t-x^\star}+
\frac{1- p_\emptyset + p_\emptyset\hat{p}}{n}\sum_{i=1}^n 
\frac{ 1-p_i}{p_i}\gamma^2  \eta_i  
\sqn{u_i^{t} -u_i^\star }\\
&\quad + (1- p_\emptyset + p_\emptyset\hat{p})\sqn{w^t-w^\star}-(1- p_\emptyset + p_\emptyset\hat{p})\gamma^2\sqn{\tilde{\nabla}\g(\hat{x}^{t})-\tilde{\nabla}\g(\hat{x}^\star)+\frac{1}{n} \sum_{j=1}^n (u_j^t-u_j^\star)}\\
&\quad-2\gamma(1- p_\emptyset + p_\emptyset\hat{p}) \langle \tilde{\nabla}\g(\hat{x}^{t})-\tilde{\nabla}\g(\hat{x}^\star), \hat{x}^{t} -x^\star\rangle +
 \frac{1- p_\emptyset + p_\emptyset\hat{p}}{n} \sum_{i=1}^n \gamma^2  \eta_i \sqn{u_{i}^t-u_i^\star}\\
&\quad-2\gamma \frac{1- p_\emptyset +p_\emptyset\hat{p}}{n}\sum_{i=1}^n   \langle y_{i}^{t+1}-x^\star,\tilde{\nabla} \hi(y_i^{t+1}) - u_i^\star\rangle.
\end{align*}
We have, since $\gamma<\frac{2}{L_\ff}$,
\begin{align*}
\sqn{w^t-w^\star}&=\sqn{x^t-x^\star}-2\gamma\langle x^t-x^\star,\nabla\ff(x^t)-\nabla\ff(x^\star)\rangle+\gamma^2\sqn{\nabla\ff(x^t)-\nabla\ff(x^\star)}\\
&\leq \sqn{x^t-x^\star}-\gamma(2-\gamma L_\ff)\langle x^t-x^\star,\nabla\ff(x^t)-\nabla\ff(x^\star)\rangle.
\end{align*}
Hence,
\begin{align*}
&\Expc{\sqn{x^{t+1}-x^\star}} +\frac{1- p_\emptyset + p_\emptyset\hat{p}}{n}\gamma^2 \sum_{i=1}^n \frac{\eta_i}{p_i}
\Expc{\sqn{u_i^{t+1}-u_i^\star}}\notag\\
&\leq \sqn{x^t-x^\star} +
\frac{1- p_\emptyset + p_\emptyset\hat{p}}{n}\gamma^2 \sum_{i=1}^n 
\frac{ \eta_i}{p_i}
\sqn{u_i^{t} -u_i^\star }\\
&\quad 
-\gamma(2-\gamma L_\ff)(1- p_\emptyset + p_\emptyset\hat{p})\langle x^t-x^\star,\nabla\ff(x^t)-\nabla\ff(x^\star)\rangle\\
&\quad-(1- p_\emptyset + p_\emptyset\hat{p})\gamma^2\sqn{\tilde{\nabla}\g(\hat{x}^{t})-\tilde{\nabla}\g(\hat{x}^\star)+\frac{1}{n} \sum_{j=1}^n (u_j^t-u_j^\star)}\\
&\quad-2\gamma(1- p_\emptyset + p_\emptyset\hat{p}) \langle \tilde{\nabla}\g(\hat{x}^{t})-\tilde{\nabla}\g(\hat{x}^\star), \hat{x}^{t} -x^\star\rangle \\
&\quad-2\gamma \frac{1- p_\emptyset +p_\emptyset\hat{p}}{n}\sum_{i=1}^n   \langle y_{i}^{t+1}-x^\star,\tilde{\nabla} \hi(y_i^{t+1}) - u_i^\star\rangle.
\end{align*}
Let $T\geq 0$. By unrolling the recursion above, we obtain
\begin{align*}
&\Exp{\sqn{x^{T+1}-x^\star}} +\frac{1- p_\emptyset + p_\emptyset\hat{p}}{n}\gamma^2 \sum_{i=1}^n \frac{\eta_i}{p_i}
\Exp{\sqn{u_i^{T+1}-u_i^\star}}\notag\\
&\quad +\gamma(2-\gamma L_\ff)(1- p_\emptyset + p_\emptyset\hat{p})\sum_{t=0}^{T}\Exp{\langle x^t-x^\star,\nabla\ff(x^t)-\nabla\ff(x^\star)\rangle}\\
&\quad+(1- p_\emptyset + p_\emptyset\hat{p})\gamma^2\sum_{t=0}^{T}\Exp{\sqn{\tilde{\nabla}\g(\hat{x}^{t})-\tilde{\nabla}\g(\hat{x}^\star)+\frac{1}{n} \sum_{j=1}^n (u_j^t-u_j^\star)}}\\
&\quad+2\gamma(1- p_\emptyset + p_\emptyset\hat{p})\sum_{t=0}^{T}\Exp{ \langle \tilde{\nabla}\g(\hat{x}^{t})-\tilde{\nabla}\g(\hat{x}^\star), \hat{x}^{t} -x^\star\rangle} \\
&\quad+2\gamma \frac{1- p_\emptyset +p_\emptyset\hat{p}}{n}\sum_{t=0}^{T}\sum_{i=1}^n  \Exp{ \langle \tilde{\nabla} \hi(y_i^{t+1}) - u_i^\star,y_{i}^{t+1}-x^\star\rangle}\\
&\leq \sqn{x^0-x^\star} +
\frac{1- p_\emptyset + p_\emptyset\hat{p}}{n}\gamma^2 \sum_{i=1}^n 
\frac{ \eta_i}{p_i}
\sqn{u_i^{0} -u_i^\star }.
\end{align*}
Thus, by choosing $t$ uniformly at random in $\{0,\ldots,T\}$, we have 
\begin{align*}
&\gamma(2-\gamma L_\ff)(1- p_\emptyset + p_\emptyset\hat{p})\Exp{\langle x^t-x^\star,\nabla\ff(x^t)-\nabla\ff(x^\star)\rangle}\\
&\quad+(1- p_\emptyset + p_\emptyset\hat{p})\gamma^2\Exp{\sqn{\tilde{\nabla}\g(\hat{x}^{t})-\tilde{\nabla}\g(\hat{x}^\star)+\frac{1}{n} \sum_{j=1}^n (u_j^t-u_j^\star)}}\\
&\quad+2\gamma(1- p_\emptyset + p_\emptyset\hat{p})\Exp{ \langle \tilde{\nabla}\g(\hat{x}^{t})-\tilde{\nabla}\g(\hat{x}^\star), \hat{x}^{t} -x^\star\rangle} \\
&\quad+2\gamma \frac{1- p_\emptyset +p_\emptyset\hat{p}}{n}\sum_{i=1}^n  \Exp{ \langle \tilde{\nabla} \hi(y_i^{t+1}) - u_i^\star,y_{i}^{t+1}-x^\star\rangle}\\
&\leq \frac{1}{T+1}\left(\sqn{x^0-x^\star} +
\frac{1- p_\emptyset + p_\emptyset\hat{p}}{n}\gamma^2 \sum_{i=1}^n 
\frac{ \eta_i}{p_i}
\sqn{u_i^{0} -u_i^\star }\right).
\end{align*}

\section{Comparison of Our Linear Rate with Existing Results for the Davis--Yin algorithm}\label{seccdy}

We consider the case $n=1$ and $|\Omega^t|\equiv 1$, so that $p_\emptyset=0$. Then \algno reverts to the \algn{Davis--Yin algorithm}, as discussed in Section \ref{secsota}.  We can compare our linear rates to known rates in the literature, keeping in mind that they are not for the same Lyapunov function.
Our rate in \eqref{eqrho6} is
\begin{align*}
\rho&=\max\left(
\frac{\max(1-\gamma\mu_{\ff},\gamma L_{\ff}-1)^2}{(1+\gamma\mu_{\g})(1+\gamma \mu_{\dr{h}})},
\frac{\gamma(L_{\dr{h}}+\mu_{\dr{h}})}{\gamma(L_{\dr{h}}+\mu_{\dr{h}})+2}
\right)
\end{align*}
(we write $\h$ instead of $\dr{h_1}$). 
\citet[Theorem 1]{yir22} (with $L_\g = +\infty$ since we do not assume $\g$ to be smooth) provide the rate
\begin{equation*}
\rho_\mathrm{DY}= \max\left(\frac{\gamma^2 \mu_{\dr{h}}^2}{(1+\gamma\mu_{\dr{h}})^2}+\frac{\max(1-\gamma\mu_{\ff},\gamma L_{\ff}-1)^2}{(1+2\gamma\mu_\g)(1+\gamma\mu_{\dr{h}})^2},
\frac{\gamma^2 L_{\dr{h}}^2}{(1+\gamma L_{\dr{h}})^2}+\frac{\max(1-\gamma\mu_{\ff},\gamma L_{\ff}-1)^2}{(1+2\gamma \mu_\g)(1+\gamma L_{\dr{h}})^2}
\right).
\end{equation*}
A comparison with our rate is not obvious. But when $\mu_{\dr{h}}=0$, Yi and Ryu showed that 
$\rho_\mathrm{DY}\leq \rho$, with 
 \begin{align*}
\rho_\mathrm{DY}&= \max\left(\frac{\max(1-\gamma\mu_{\ff},\gamma L_{\ff}-1)^2}{1+2\gamma\mu_\g},
\frac{\gamma^2 L_{\dr{h}}^2}{(1+\gamma L_{\dr{h}})^2}+\frac{\max(1-\gamma\mu_{\ff},\gamma L_{\ff}-1)^2}{(1+2\gamma \mu_\g)(1+\gamma L_{\dr{h}})^2}
\right)\\
&\geq \max\left(\frac{\max(1-\gamma\mu_{\ff},\gamma L_{\ff}-1)^2}{1+2\gamma\mu_\g},
\frac{\gamma L_{\dr{h}}}{ \gamma L_{\dr{h}}+2+\frac{1}{\gamma L_{\dr{h}}}}
\right),
\end{align*}
whereas we have
 \begin{equation*}
\rho=   \max\left( \frac{\max(1-\gamma\mu_{\ff},\gamma L_{\ff}-1)^2}{1+\gamma\mu_{\g}},
 \frac{\gamma L_{\dr{h}} }{\gamma L_{\dr{h}}+2 }\right).
\end{equation*}
We see that the improvement is essentially the factor 2 in the term $1+2\gamma\mu_\g$ in $\rho_\mathrm{DY}$ versus $1+\gamma\mu_\g$ in $\rho$.
 
 When $\ff=0$, the \algn{Davis--Yin algorithm} becomes the \algn{Douglas--Rachford algorithm}. 
 In the case $\mu_\g=0$, a tight rate 
 for this  algorithm 
 is proved in \citet[Theorem 2]{gis17}:
\begin{align*}
\rho_\mathrm{DR}&=\max\left(\frac{1}{(\gamma\mu_\h+1)^2},\frac{\gamma^2 L_{\h}^2}{(\gamma L_\h+1)^2}\right).
\end{align*}
Our rate is
\begin{equation*}
\rho=   \max\left(
\frac{L_{\dr{h}}+\mu_{\dr{h}}}{L_{\h}+\mu_{\dr{h}} + 2\gamma L_{\dr{h}}\mu_{\dr{h}} },
\frac{\gamma (L_{\h}+\mu_\h) }{\gamma (L_{\h}+\mu_\h)+2 }\right).
\end{equation*}
Comparing the first factor of $\rho$ and $\rho_\mathrm{DR}$, we have
\begin{align*}
\frac{1}{(1+\gamma\mu_\h)^2}\leq \frac{1}{1+2\gamma \mu_\h}\leq \frac{L_{\h}+\mu_{\h}}{L_{\h}+\mu_{\h} + 2\gamma L_{\h}\mu_{\h} }\leq \frac{1}{1+\gamma \mu_\h},
\end{align*}
and in practical cases with $\gamma\mu_\h\ll 1$ and $\mu_\h\ll L_\h$, both terms are very close to $\frac{1}{1+2\gamma \mu_\h}$. Comparing the second factor of $\rho$ and $\rho_\mathrm{DR}$, we have
\begin{align*}
\frac{\gamma^2 L_{\h}^2}{(\gamma L_\h+1)^2}=\frac{\gamma L_{\h}}{\gamma L_\h+2+\frac{1}{\gamma L_\h}}\leq
\frac{\gamma L_{\h} }{\gamma L_{\h}+2 } \leq  \frac{\gamma (L_{\h}+\mu_\h) }{\gamma (L_{\h}+\mu_\h)+2 } \leq  \frac{\gamma L_{\h} }{\gamma L_{\h}+1 } ,
\end{align*}
and in practical cases with $\gamma L_\h\gg 1$ and $\mu_\h\ll L_\h$, both terms are very close to $\frac{\gamma L_{\h} }{\gamma L_{\h}+2 }$. 
Thus, our rate $\rho$ is very close to the tight rate  $\rho_\mathrm{DR}$.

\end{document}